\documentclass[a4paper,12pt,twoside, reqno]{amsart}

\usepackage{enumitem}
\usepackage{amsmath}
\usepackage{amssymb}
\usepackage{amsfonts}
\usepackage{amscd}
\usepackage{hyperref}

\newtheorem{theorem}{Theorem}
\newtheorem{lemma}{Lemma}

\newtheorem{fact}{Fact}
\newtheorem{corollary}{Corollary}
\theoremstyle{definition}
\newtheorem{definition}{Definition}
\newtheorem{remark}{Remark}

\def\Z{{\mathbb Z}}
\def\N{{\mathbb N}}
\def\Q{{\mathbb Q}}
\def\MA{{\mathbb A}}
\def\MB{{\mathbb B}}
\def\MC{{\mathbb C}}

\newcommand{\nsZ}{\widetilde{\mathbb Z}}
\newcommand{\nsA}{\widetilde{\mathbb A}}
\newcommand{\nsB}{\widetilde{\mathbb B}}

\newcommand{\ncl}{\mathrm{ncl}}
\newcommand{\Aut}{\mathrm{Aut}}
\newcommand{\Th}{\mathrm{Th}}

\title[Baumslag\,--\,Solitar groups and regular bi-interpretability]{Groups elementarily equivalent to metabelian Baumslag\,--\,Solitar groups and regular bi-interpretability}
\author{E.\,Daniyarova, A.\,Myasnikov}

\begin{document}

\begin{abstract} 
We prove that metabelian Baumslag\,--\,Solitar group $BS(1,k)$, $k>1$, is (strongly) regularly bi-interpretable with the ring of integers $\Z$, and describe in algebraic terms all groups that are elementarily equivalent to $BS(1,k)$.
\end{abstract}

\maketitle


\tableofcontents

\section{Introduction}

In this paper we describe all groups elementarily equivalent to a given metabelian Baumslag\,--\,Solitar group $BS(1,k)$, $k>1$, in the standard group theory language $\mathcal L_{\mathrm Gr} = \{\cdot, \,^{-1}, e\}$. It turns out these groups have a very particular algebraic structure, they come, akin to the classical situation in arithmetic, as the non-standard models of the initial group $BS(1,k)$. Note that many groups $G$ (for example, all countable groups with arithmetical multiplication table) have non-standard models $G(\tilde \Z)$, for each ring  of non-standard integers $\tilde \Z \equiv \Z$ (in particular, $G = G(\Z)$), though in general the elementary theory $\Th(G)$ may have some other models besides the non-standard ones. The important property of the non-standard models $G(\tilde \Z)$ is that they are  equivalent to $G$ (and to each other) in the so called \emph{non-standard weak second order logic}, which is much stronger then the first-order one. Hence they more closely resemble the initial group $G$ than the other models of $\Th(G)$.
In our  case when $G = BS(1,k)$, all models of $\Th(G)$ are precisely the non-standard ones. To prove this, we developed a new approach to the first-order classification problem  based on recent advances of theory of interpretability, especially on the techniques of bi-interpretability (see, for example,~\cite{AKNS,Khelif,Nies1,Nies2,KhMS,Th_int2}). Our main technical result is that every group $BS(1,k)$ with $k >1$ is regularly bi-interpretable with $\Z$. The notion of regular interpretability is somewhat in between the absolute interpretability and interpretability with parameters. The former one fits extremely well for questions related to the first-order equivalence, but occurs rather rarely. The interpretability with parameters is a more common one, but the parameters are usually not in the language of $G$, which presents a problem when describing models of $\Th(G)$. The regular interpretability (bi-interpretability) combines the useful properties of both: it can be used in first-order classification problem in a way similar to the absolute one, and occurs almost as often as interpretability with parameters. However, proving regular interpretability (bi-interpretability) is more demanding, it requires one to go deeper into the model-theoretic properties of the structure $G$.  

The logical and algorithmic properties of the groups $BS(1,k)$ were thoroughly studied. It does not happen by chance; on the one hand, these groups provide one of the typical fundamental examples of finitely generated metabelian groups; on the other hand, quite often they play the role of principal obstacles in a study of model-theoretic, algebraic, and algorithmic properties of other types of groups. Noskov showed in~\cite{Noskov} that the elementary theories $\Th(G)$ of finitely generated solvable not virtually abelian groups $G$ in the language $\mathcal L_{\mathrm Gr}$ are undecidable, in particular,  $\Th(BS(1,k))$ is undecidable. In fact, he proved in this paper, that the ring $\Z$ is interpretable (with parameters) in $BS(1,k)$. In~\cite{Khelif} Kh\'elif outlined an argument that shows that $BS(1,k)$ and $\Z$ are bi-interpretable (with parameters) in each other. This implies that $BS(1,k)$ is a prime model of its theory $\Th(BS(1,k))$ and it is QFA, that is, any other finitely generated group $H$ such that $H \equiv BS(1,k)$ is, in fact, isomorphic to $BS(1,k)$ (see papers~\cite{Nies1,Nies2} by Nies for this). In~\cite{KhMS} Kharlampovich, Myasnikov and Sohrabi introduced a notion of a \emph{rich structure} and proved that the group $BS(1,k)$ is rich, i.\,e., every formula of the weak second order logic in the language $\mathcal L_{Gr}$ is equivalent in $BS(1,k)$ to some first-order logic formula in the language $\mathcal L_{\mathrm Gr}$. This implies that many interesting group-theoretic properties of $BS(1,k)$ are definable by first-order formulas in $BS(1,k)$ (see~\cite{KhMS}). In contrast to the first-order theory, the universal theory of $BS(1,k)$  in the language $\mathcal L_{\mathrm Gr}$ is decidable, see~\cite{GupTim} by Gupta and Timoshenko. In~\cite{Rom} Romanovskii described in algebraic terms all finitely generated groups universally equivalent to $BS(1,k)$. Note that again, unlike first-order equivalence, there are infinitely many finitely generated metabelian groups universally equivalent to $BS(1,k)$. Elementary equivalence of arbitrary (not necessarily metabelian) Baumslag\,--\,Solitar groups was studied in~\cite{CK}, where Casals-Riuz and Kazachkov showed that in this class of groups elementary equivalence implies isomorphism.

Our solution of the first-order classification problem for metabelain Baumslag\,--\,Solitar groups $BS(1,k)$, $k >1,$ is based on the following important technical result of Section~\ref{subsec:4.5}.

\medskip\noindent 
{\bf Theorem~\ref{prop2}.} 
{\it If $k>1$, then the group $BS(1,k)$ is regularly strongly bi-interpretable with $\Z$.} 

\medskip
The notion of regular bi-interpretation is rather technical, we explain it in detail in Section~\ref{se:2}, which also contains all the necessary notation and terminology. To prove Theorem~\ref{prop2} we need various facts from number theory and commutative algebra on the ring of Laurent polynomials 
$\Z[x,x^{-1}]$ and the ring $\Z[1/k]$, which we prove in Section~\ref{se:3}. Using these results, we study in detail in Section~\ref{subsec:4.1} various sets and predicates definable in $BS(1,k)$, which are crucial for the model theory of $BS(1,k)$. Note that throughout the paper, we treat the groups $BS(1,k)$ as semi-direct products $\Z[1/k]\rtimes \Z$, where $\Z[1/k]$ and $\Z$ are viewed as the additive groups of the corresponding rings. 

Theorem~\ref{prop2} implies that all groups elementarily equivalent to $BS(1,k)$, $k>1$, are precisely the non-standard models of $BS(1,k)$. For this we need some recent results on regular bi-interpretations, which we outline in Section~\ref{subsec:2.3}.  

In Section~\ref{subsec:5.2} we introduce and study the algebraic structure of \emph{non-standard models} of the group $BS(1,k)$. We start with an absolute interpretation $BS(1,k) \simeq \Delta(\Z)$ of $BS(1,k)$ in $\Z$. By definition, the non-standard version of $BS(1,k)$ is the group $\Delta(\nsZ)$ obtained by the interpretation $\Delta$ applied to a ring $\nsZ$ with $\nsZ \equiv \Z$. It turns out that the group $\Delta(\nsZ)$ does not depend on the choice of interpretation $\Delta$, it depends only on $\nsZ$. In fact, for all such interpretations $\Delta$, the group $\Delta(\nsZ)$ is isomorphic to the group $BS(1,k,\nsZ)$, which is a semi-direct product 
$$
BS(1,k,\nsZ) = \nsZ[1/k^{\nsZ}]\rtimes \nsZ,
$$
where $\nsZ[1/k^{\nsZ}]$ and $\nsZ$ are viewed as the additive groups of the corresponding rings. We need to describe algebraically the ring $\nsZ[1/k^{\nsZ}]$. The standard exponentiation $x=y^z$, where $x,y, z\in \Z$, $z \geq 0$,  is definable in $\Z$, hence the same formula defines a non-standard exponentiation in $\nsZ$, which satisfies the standard high school identities. In particular, for any $k\in \N$ and $i\in \nsZ^+ =\{u \in \nsZ \mid u >  0\}$ there exists a unique element $k^i\in \nsZ^+$, such that the classical identities $k^{i_1}\cdot k^{i_2}=k^{i_1+i_2}$, $i_1,i_2\in\nsZ^+$ hold.
 Thus, the set $S=\{k^i,i\in \nsZ^+\}\cup\{ 1\}$ is a multiplicative subset of $\nsZ$, so one can form  the ring of fractions $S^{-1}\nsZ$, which we denote by $\nsZ[1/k^{\nsZ}]$. In $\nsZ[1/k^{\nsZ}]$ the elements $k^i, i \in \nsZ^+$ are invertible, and we may write $k^{-i}=1/k^i$. Now, any element $x\in \nsZ[1/k^{\nsZ}]$ can be written as $x=z\cdot k^i$, where $z,i\in \nsZ$.  In the semi-direct product $\nsZ[1/k^{\nsZ}]\rtimes \nsZ$ above,  $\nsZ$ acts on $\nsZ[1/k^{\nsZ}]$ via a homomorphism  $\varphi\colon \nsZ \to \Aut(\nsZ[1/k^{\nsZ}])$, where  for $m \in \nsZ$ the automorphism $\varphi(m) \in \Aut(\nsZ[1/k^{\nsZ}])$ is defined as
$$
\varphi(m)\colon x \to  x\cdot k^{-m}, \ \ x\in\nsZ[1/k^{\nsZ}].
$$

We believe that the groups $BS(1,k,\nsZ)$ are interesting in their own right from algebraic, geometric, and even algorithmic (via Turing machines over $\nsZ$) view-points. We leave it for future research; however, we describe here one of their interesting properties. 
In~\cite{Lyn} Lyndon introduced the notion of a group with an exponentiation in an associative unitary ring $A$. Later, Myasnikov and Remeslennikov added one more axiom to Lyndon's definition, so the refined $A$-groups can be viewed as non-commutative $A$-modules~\cite{MR1} (termed $A$-groups). Recall that an $A$-group $G$ is a group that comes equipped with an exponentiation $g \to g^a$, where $g \in G, a \in A$ which satisfies the following axioms:

\begin{enumerate}
\item    $\forall\,g \; (g^1=g\:\wedge \:g^0=e)$; $e^\alpha=e$, $\alpha \in A$;
\item $\forall\,g \;  (g^{\alpha+\beta}=g^{\alpha}g^{\beta} \:\wedge\: g^{\alpha\beta}=(g^{\alpha})^{\beta})$, $\alpha,\beta \in A$;
\item   $\forall\,g \:\forall\,h\; ((h^{-1}gh)^\alpha=h^{-1}g^{\alpha}h)$, $\alpha \in A$;
\item 
    $\text{$(M\!R$-axiom)}\quad \forall\,g\:\forall\,h \; (gh=hg\;\longrightarrow \;(gh)^{\alpha}=g^{\alpha}h^{\alpha})$, $\alpha \in A$.
\end{enumerate}

\medskip\noindent
{\bf Theorem~\ref{th:exponent}.} 
{\it 
    For every $k>1$ and $\nsZ \equiv \Z$ the group $BS(1,k,\nsZ)$ is a $\nsZ$-group. 
}

\medskip

In~\cite{MR2} a general notion of a tensor $A$-completion $G \otimes A$ of a group $G$ by a ring $A$ was introduced and its algebraic properties for various groups $G$ were studied. In~\cite{AMN} authors went further and defined a tensor product $G \otimes_{\mathcal V} A$ in a variety  of groups $\mathcal V$. In particular, tensor completion $G \otimes_{\mathcal M_2} A$  is defined in the variety of metabelian groups $\mathcal M_2$. It will be very interesting to compare $\nsZ$-groups $BS(1,k,\nsZ)$ and $BS(1,k) \otimes_{\mathcal M_2} \nsZ$.  

Finally, in Section~\ref{subsec:5.3} we prove the following result that describes all groups elementarily equivalent to $BS(1,k)$.

\medskip\noindent
{\bf Theorem~\ref{th:FO classification}.} 
{\it 
Let $k$ be an integer, $k>1$. Then for any ring $\nsZ\equiv \Z$ one has 
$$
BS(1,k)\equiv BS(1,k,\nsZ)
$$ 
and, conversely, every group  elementarily equivalent to $BS(1,k)$ has the form $BS(1,k,\nsZ)$ for some $\nsZ\equiv \Z$. Moreover, this description of models of $\Th(BS(1,k))$ is a bijection: for any $\nsZ_1\equiv\Z\equiv\nsZ_2$ one has 
$$
BS(1,k,\nsZ_1) \simeq BS(1,k,\nsZ_2) \Longleftrightarrow \nsZ_1\simeq \nsZ_2.
$$
}

\section{Interpretability and non-standard models}
\label{se:2}

In this section, we discuss various types of interpretability of structures; some of them are well known, others are new. Each of them serves different purposes. Our focus is mostly on the regular interpretability, which suits well to the classical first-order classification problem.

In the following we work only with groups and rings, so to simplify arguments, we consider only algebraic structures without predicates in the signature. Even though all the model-theoretic results below hold for arbitrary languages.

Let $L$ be a functional language with a set of functional (operational) symbols $\{f, \ldots\}$ together with their arities $n_f \in \N$, and a set of constant symbols $\{c, \ldots\}$. We write $f(x_1, \ldots,x_n)$ to show that $n_f = n$. The standard language of groups $\{\cdot\,,\,^{-1},e\}$ includes the symbol $\cdot$ for binary operation of multiplication, the symbol $^{-1}$ for unary operation of inversion, and the symbol $e$ for the group identity; the standard language of unitary rings is $\{+,\,\cdot\,,0,1\}$. An interpretation of a constant symbol $c$ in a set $A$ is an element $c^A\in A$. For a functional symbol $f$ an interpretation in $A$ is a function $f^A\colon A^{n_f}\to A$. 

An algebraic structure in the language $L$ (an $L$-structure) with the base set $A$ is denoted by $\MA = \langle A; L\rangle$, or simply by $\MA = \langle A; f, \ldots,c, \ldots \rangle$, or by $\MA = \langle A; f^A, \ldots,c^A, \ldots \rangle$. For a given structure $\MA$ by $L(\MA)$ we denote the language of $\MA$.   

We usually denote variables by small letters $x,y,z, a,b, u,v, \ldots$, while the same symbols with bars $\bar x, \bar y,  \ldots$ denote tuples of the corresponding variables, say $\bar x = (x_1, \ldots,x_n)$, and furthermore, $\bar{\bar x}$ is a tuple of tuples $\bar{\bar x}=(\bar x_1,\ldots,\bar x_m)$, where $|\bar x_i|=n$.

\subsection{Definable sets and interpretability}

Let $\MB = \langle B; L(\MB)\rangle$ be an algebraic structure. A subset $X \subseteq B^n$ is called {\em $0$-definable} (or {\em absolutely definable}, or {\em definable without parameters}) in $\MB$ if there is a first-order formula $\phi(x_1,\ldots,x_n)$ in the language $L(\MB)$ such that
$$
X = \{(b_1,\ldots,b_n) \in B^n \mid \MB \models \phi(b_1, \ldots,b_n)\}.
$$ 
The set $X$ is denoted by $\phi(\MB)$. Fix a set of elements $P\subseteq B$. If $\psi(x_1,\ldots,x_n,y_1,\ldots,y_k)$ is a formula in $L(\MB)$ and $\bar{p}=(p_1,\ldots,p_k)$ is a tuple of elements from $P$, then the set 
$$
\{(b_1,\ldots,b_n) \in B^n \mid \MB \models \psi(b_1,\ldots,b_n, p_1,\ldots,p_k)\}
$$ 
is called \emph{definable} in $\MB$ (\emph{with parameters} $\bar{p}$). It is denoted by $\psi(\MB,\bar p)$ or $\psi(B^n,\bar p)$.

Equivalently, $X\subseteq B^n$ is definable in $\MB$ if and only if there exists a finite set of parameters $P\subseteq B$ and a formula $\psi(x_1,\ldots,x_n)$ in the language $L(\MB)\cup P$ such that $X=\psi(\MB_P)$, here the algebraic structure $\MB_P=\langle B;L(\MB)\cup P\rangle$ is obtained from $\MB$ by adding new constant symbols from $P$ into the language. 

Let $X\subseteq B^n$, $Y\subseteq B^m$ be sets and $\sim_X$, $\sim_Y$ be equivalence relations on $X$, $Y$, respectively. A map $g\colon X/{\sim_X}\to Y/{\sim_Y}$ is called {\em definable} in $\MB$, if the following set (termed the {\em preimage in $\MB$ of the graph of $g$})
\begin{multline*}
\{(b_1,\ldots,b_n,c_1,\ldots,c_m) \in B^{n+m}\mid  \\ 
g((b_1,\ldots,b_n)/{\sim_X})=(c_1,\ldots,c_m)/{
\sim_Y},\\
(b_1,\ldots,b_n)\in X, (c_1,\ldots,c_m)\in Y\}
\end{multline*}
is definable in $\MB$. Note that $g$ may  also be a constant (i.\,e.,  a function of arity $0$).

Note that definable sets in arithmetic $\N$, the so-called arithmetic sets, are well-studied. In particular, it is known that every computably enumerable set in $\N$ is definable. The same holds for the ring of integers $\Z$. Furthermore, every element of $\N$ (or $\Z$) is absolutely definable in 
$\N$ ($\Z$), so every definable set in $\N$ ($\Z$) is absolutely definable.

\begin{definition}  \label{de:interpretation}
An algebraic  structure $\MA = \langle A;f,\ldots,c,\ldots\rangle$ is $0$-{\em interpretable} (or {\em absolutely interpretable}, or {\em interpretable without parameters}) in an algebraic structure $\MB=\langle B;L(\MB)\rangle$ if the following conditions hold:
\begin{enumerate}[label=\arabic*)]
\item there is a subset $A^\ast  \subseteq B^n$  $0$-definable in $\MB$, 
\item there is an equivalence relation $\sim$ on $A^\ast$ $0$-definable in $\MB$, 
\item there are interpretations $f^A, \ldots$, $c^A, \ldots$ of the symbols $f, \dots, c, \ldots$ on the quotient set $A^\ast /{\sim}$, all $0$-definable in $\MB$, 
\item the structure $\MA^\ast=\langle A^\ast /{\sim}; f^A, \ldots,c^A, \ldots \rangle$ is $L(\MA)$-isomorphic to $\MA$.
\end{enumerate}
\end{definition}

A more general case is the interpretation with parameters.

\begin{definition} \label{de:interpretation_P}
An algebraic structure $\MA=\langle A;f,\ldots,c,\ldots\rangle$  is {\em interpretable} (with parameters) in an algebraic structure $\MB=\langle B;L(\MB)\rangle$ if there is a finite subset of parameters $P\subseteq B$, such that $\MA$ is absolutely interpretable in $\MB_P$. In this case, we write $\MA\rightsquigarrow\MB$. 
\end{definition}

The structure $\MA^\ast$ from Definitions~\ref{de:interpretation}, \ref{de:interpretation_P} is completely described by the tuple of parameters $\bar p$ (all the parameters used in the formulas from 1), 2), and 3) above) and the following set of formulas in the language $L(\MB)$:
\begin{equation} \label{eq:code}
\Gamma =  \{U_\Gamma(\bar x,\bar y), E_\Gamma(\bar x, \bar x^\prime,\bar y), Q_\Gamma(\bar x_1, \ldots,\bar x_{t_Q},\bar y) \mid Q \in L(\MA)\},
\end{equation}
 where $\bar x$, $\bar x^\prime$ and $\bar x_i$ are $n$-tuples of variables and $\bar y = (y_1, \ldots,y_k)$ is a tuple of extra variables for parameters $\bar p$. Namely,  $U_\Gamma$  defines in $\MB$ a set $A_\Gamma  = U_\Gamma(B^n,\bar p)  \subseteq B^n$ (the set $A^\ast$ in Definition~\ref{de:interpretation}), $E_\Gamma$ defines an equivalence relation $\sim_\Gamma$ on $A_\Gamma$ (the equivalence relation $\sim$ in Definition~\ref{de:interpretation}), and the formulas $Q_\Gamma$ define the preimages in $\MB$ of the graphs for constants and functions $Q\in L(\MA)$ on the quotient set $A_\Gamma/{\sim_\Gamma}$ in such a way that the structure $\Gamma(\MB,\bar p) = \langle A_\Gamma/{\sim_\Gamma}; L(\MA) \rangle $ is isomorphic to $\MA$. Note that we interpret a constant $c \in L(\MA)$ in the structure $\Gamma(\MB,\bar p)$ by the $\sim_\Gamma$-equivalence  class of some tuple $\bar b_c \in A_\Gamma$ defined in $\MB$ by the formula $c_\Gamma(\bar x,\bar p)$. The number $n$ is called the {\em dimension} of $\Gamma$, denoted $n = \dim \Gamma$. We refer to $k$ as the {\em parameter dimension} of $\Gamma$ and denote it by $\dim_{par}\Gamma$. 

We refer to $\Gamma$  as the {\em interpretation code} or just the {\em code} of the interpretation $\MA\rightsquigarrow\MB$. Sometimes we identify the interpretation $\MA\rightsquigarrow\MB$ with its code $\Gamma$ or $(\Gamma,\bar p)$ (if the interpretation has parameters $\bar p$). We may also write $\MA\stackrel{\Gamma}{\rightsquigarrow}\MB$, $\MA\stackrel{\Gamma, \bar p}{\rightsquigarrow}\MB$, $\MA\simeq\Gamma(\MB)$ or $\MA \simeq \Gamma(\MB,\bar p)$. To stress that the interpretation $\Gamma$ is absolute we write $\MA\simeq\Gamma(\MB,\emptyset)$. 

By $\mu_\Gamma$ we denote a surjective map $A_\Gamma \to A$ that gives rise to an isomorphism $\bar \mu_\Gamma\colon \Gamma(\MB,\bar p) \to \MA$. We refer to $\mu_\Gamma$ as the {\em coordinate map} of the interpretation $(\Gamma,\bar p)$. Note that such $\mu_\Gamma$ may not be unique because $\bar\mu_\Gamma$ is defined up to an automorphism of $\MA$. For this reason, the coordinate map $\mu_\Gamma$ is sometimes added to the interpretation notation: $(\Gamma,\bar p, \mu_\Gamma)$.  A coordinate map $\mu_\Gamma\colon A_\Gamma \to A$ gives rise to a map $\mu_\Gamma^m\colon A^m_\Gamma\to A^m$ on the corresponding Cartesian powers, which we often denote by  $\mu_\Gamma$.

When $\MA$ is definable in $\MB$, so the formula $E_\Gamma$ defines the identity relation $(x_1=x^\prime_1)\wedge\ldots\wedge(x_n=x^\prime_n)$, the surjection $\mu_\Gamma$ is injective. In this case, $(\Gamma,\bar p, \mu_\Gamma)$ is also called an {\em injective interpretation}.

\begin{definition}\label{def:regular_int}
We say that $\MA$ is {\em regularly interpretable} in $\MB$ if there exists a code $\Gamma$~\eqref{eq:code} and an $L(\MB)$-formula $\phi(y_1,\ldots,y_k)$, such that $\phi(\MB)\ne\emptyset$ and for each $\bar p=(p_1,\ldots,p_k)\in \phi(\MB)$ $\Gamma(\MB,\bar p)$ gives an interpretation $\MA\simeq\Gamma(\MB,\bar p)$ of $\MA$ in $\MB$.
\end{definition}

We denote by $(\Gamma,\phi)$  the regular interpretation with the code $\Gamma$ and the formula $\phi$ for the parameters. If for all $\bar p\in \phi(\MB)$ the interpretation $(\Gamma,\bar p)$ is injective, then the regular interpretation $(\Gamma,\phi)$ is called {\em injective}. We write  $\MA\simeq \Gamma(\MB,\phi)$ if $\MA$ is regularly interpreted in $\MB$ by a code $\Gamma$ and a formula $\phi$.

Note that every absolute interpretation $\Gamma$ is a particular case of a regular interpretation. Indeed, one can add some fictitious variable $y$ in all the formulas of the code $\Gamma$ and put $\phi(y)=(y=y)$.

Let $\Gamma$ be the code described in~\eqref{eq:code}. The \emph{admissibility conditions $\mathcal{AC}_\Gamma(\bar y)$} (see~\cite{Th_int1}) for the code $\Gamma$ is a set of formulas in the language $L(\MB)$ with free variables $\bar y=(y_1, \ldots,y_k)$, where $k=\dim_{par}\Gamma$, such that for an arbitrary $L(\MB)$-structure $\nsB$ and an arbitrary tuple $\bar q$ in $\nsB$ the construction of $\Gamma(\nsB,\bar q)$ described in Definitions~\ref{de:interpretation}, \ref{de:interpretation_P} indeed gives an $L(\MA)$-structure if and only if $\nsB\models \mathcal{AC}_\Gamma(\bar q)$. Observe that the admissibility conditions $\mathcal{AC}_\Gamma(\bar q)$ state that the $L(\MA)$-structure $\Gamma(\nsB,\bar q)$ is well-defined, but they do not claim that $\MA\simeq\Gamma(\nsB,\bar q)$.

Suppose that $\phi(\bar y)$ is an $L(\MB)$-formula and $\nsB$ is an $L(\MB)$-structure as before. We say that an algebraic structure $\Gamma(\nsB,\phi)$ is {\em well-defined} if $\phi(\nsB)\ne \emptyset$ and for every $\bar q\in \phi(\nsB)$  $\Gamma(\nsB,\bar q)$ is an $L(\MA)$-structure and all these structures $\Gamma(\nsB,\bar q)$, $\bar q\in \phi(\nsB)$, are isomorphic to each other (again, they may not be isomorphic to $\MA$).

\subsection{The composition of interpretations}

It is known that the relation $\MA\rightsquigarrow\MB$  is transitive on algebraic structures (see, for example,~\cite{Th_int1, Hodges}). The proof of this fact is based on the notion of $\Gamma$-translation and composition of codes, which we present now. 

Let 
$$
\Gamma =  \{U_\Gamma(\bar x,\bar y), E_\Gamma(\bar x, \bar x^\prime,\bar y), Q_\Gamma(\bar x_1, \ldots,\bar x_{t_Q},\bar y) \mid Q \in L(\MA)\}
$$
be a code as above, consisting of $L(\MB)$-formulas. Then for    any formula $\varphi(x_1,\ldots,x_m)$ in the language $L(\MA)$ there is  a formula $\varphi_\Gamma(\bar x_1,\ldots,\bar x_m, \bar y)$ in the language $L(\MB)$  such that if $\MA\simeq\Gamma(\MB,\bar p)$, then for any coordinate map $\mu_\Gamma\colon A_\Gamma\to A$ one has
$$
\MA\models \varphi(a_1,\ldots, a_m) \iff \MB\models\varphi_\Gamma(\mu_\Gamma^{-1}(a_1),\ldots,\mu^{-1}_\Gamma(a_m),\bar p)
$$
for any elements $a_i\in A$ (see~\cite{Th_int1}). Here $\mu_\Gamma^{-1}(a_i)$ means an arbitrary preimage of $a_i$ under $\mu_\Gamma$.   Furthermore,  for any elements $\bar b_i\in B^n$ if $\MB\models\varphi_\Gamma(\bar b_1,\ldots,\bar b_m,\bar p)$ then  $\bar b_i\in \mu_\Gamma^{-1}(a_i)$ for some $a_i\in A$ with $\MA\models \varphi(a_1,\ldots, a_m)$.

Note that if the language $L(\MA)$ is computable then there is an algorithm that given $\phi$ computes $\phi_\Gamma$.

\begin{definition} \label{def:code-composition}
Let $L(\MA)$, $L(\MB)$, $L(\MC)$ be languages. Consider codes 
$$\Gamma =  \{U_\Gamma(\bar x,\bar y), E_\Gamma(\bar x, \bar x^\prime,\bar y), Q_\Gamma(\bar x_1, \ldots,\bar x_{t_Q},\bar y) \mid Q \in L(\MA)\}$$
as above and 
$$
\Delta=\{U_\Delta(\bar u,\bar z), E_\Delta(\bar u, \bar u^\prime,\bar z), Q_\Delta(\bar u_1, \ldots,\bar u_{t_Q},\bar z) \mid Q \in L(\MB)\}
$$
which consists of $L(\MC)$-formulas, where $|\bar u|=|\bar u^\prime|=|\bar u_i|=\dim \Delta$, $|\bar z|=\dim_{par}\Delta$. Then the {\em composition} of the codes $\Gamma$ and $\Delta$ is the code
$$
\Gamma \circ \Delta = \{U_{\Gamma\circ\Delta}, E_{\Gamma\circ\Delta}, Q_{\Gamma\circ\Delta} \mid Q \in L(\MA)\}=\{(U_\Gamma)_\Delta, (E_\Gamma)_\Delta, (Q_\Gamma)_\Delta \mid Q \in L(\MA)\},
$$
where $\dim \Gamma\circ \Delta =\dim \Gamma\cdot\dim \Delta$ and $\dim_{par}\Gamma\circ\Delta = \dim_{par}\Gamma\cdot \dim\Delta+\dim_{par}\Delta$.
\end{definition}

The following is an important technical result on the transitivity of interpretations.

\begin{lemma}[\cite{Th_int1}]\label{le:int-transitivity}
Let $\MA=\langle A;L(\MA)\rangle, \MB=\langle B;L(\MB)\rangle$ and $\MC=\langle C;L(\MC)\rangle$ be algebraic structures and $\Gamma,\Delta$ be codes as above. If $\MA\stackrel{\Gamma}{\rightsquigarrow}\MB$ and 
$\MB\stackrel{\Delta}{\rightsquigarrow}\MC$ then $\MA\stackrel{\Gamma\circ\Delta}{\rightsquigarrow}\MC$.

Furthermore, the following conditions hold:
\begin{enumerate}[label=\arabic*)]
    \item If $\bar p,\bar q$ are parameters and $\mu_\Gamma, \mu_\Delta$ are  coordinate maps of interpretations $\Gamma,\Delta$ then $(\bar{\bar p},\bar q)$, where $\bar{\bar p} \in \mu_\Delta^{-1}(\bar p)$, are parameters for $\Gamma\circ\Delta$; 
    \item\label{le:int2} For any coordinate map $\mu_\Delta$ of the interpretation $\MB\simeq\Delta(\MC,\bar q)$ and any tuple $\bar{\bar p} \in \mu_\Delta^{-1}(\bar p)$ the $L(\MA)$-structure $\Gamma\circ\Delta(\MC,(\bar{\bar p},\bar q))$ is well-defined and isomorphic to $\MA$;
    \item The $L(\MA)$-structure $\Gamma\circ\Delta (\MC, (\bar{\bar p},\bar q))$ does not depend on the choice of $\bar{\bar p} \in \mu_\Delta^{-1}(\bar p)$, when $\mu_\Delta$ is fixed; 
    \item $\mu_\Gamma\circ\mu_\Delta=\mu_{\Gamma}\circ\mu^n_{\Delta}\big|_{U_{\Gamma\circ\Delta}(\MC, (\bar{\bar p},\bar q))}$ is a coordinate map of the interpretation $\MA\simeq {\Gamma\circ\Delta}(\MC,(\bar{\bar p},\bar q))$ and any coordinate map $\mu_{\Gamma\circ\Delta}\colon U_{\Gamma\circ\Delta}(\MC, (\bar{\bar p},\bar q))\to A$ has a form $\mu_{\Gamma1}\circ\mu_\Delta$ for a suitable coordinate map $\mu_{\Gamma1}$ of the interpretation $\MA\simeq\Gamma(\MB,\bar p)$, provided $\mu_\Delta$ is fixed;
    \item If $\Gamma,\Delta$ are absolute, then $\Gamma\circ\Delta$ is absolute too; 
    \item If $\Gamma,\Delta$ are injective, then $\Gamma\circ\Delta$ is injective too; 
    \item\label{le:int7} If $\Gamma,\Delta$ are regular with the corresponding formulas $\varphi,\psi$, then $\Gamma\circ\Delta$ is regular too with the formula $\varphi_\Delta\wedge\psi$. 
\end{enumerate}
\end{lemma}

\begin{remark}
    Note that by construction the structure $\Gamma\circ\Delta (\MC,  (\bar{\bar p},\bar q))$  depends on the choice of $\mu_\Delta$, i.e., the set  $U_{\Gamma\circ\Delta}(\MC,(\bar{\bar p},\bar q))\subseteq C^{\dim \Gamma\circ\Delta}$, the relation $\sim_{\Gamma\circ\Delta}$,  the constants $c_{\Gamma\circ\Delta}$, and the functions   $f_{\Gamma\circ\Delta}$ in general depend on $\mu_\Delta$. On the other hand, the structure $\Gamma\circ\Delta (\MC,  (\bar{\bar p},\bar q))$ does not depend on the choice of $\mu_\Gamma$.
\end{remark}

\begin{remark}\label{remark2}
Let $(\Gamma,\phi), (\Delta,\psi)$ be regular interpretations as in~\ref{le:int7}. Then for any tuples of parameters $\bar p\in \varphi(\MB)$ and $\bar q\in \psi(\MC)$ and any coordinate map $\mu_\Delta$ of the interpretation $\MB\simeq\Delta(\MC,\bar q)$, and any preimages $\bar{\bar p} \in \mu_\Delta^{-1}(\bar p)$ the tuple $(\bar {\bar p},\bar q)$ is in $\varphi_\Delta\wedge\psi(\MC)$. Also the $L(\MA)$-structure $\Gamma\circ\Delta(\MC,(\bar {\bar p},\bar q))$ is well-defined and isomorphic to $\MA$, by item~\ref{le:int2}. Conversely, any tuple $\bar r\in \varphi_\Delta\wedge\psi(\MC)$ has a form $\bar r=(\bar{\bar p},\bar q)$, where $\bar q\in \psi(\MC)$, $\bar{\bar p}\in \varphi_\Delta(\MC,\bar q)$. For any coordinate map $\mu_\Delta$ of the interpretation $\MB\simeq\Delta(\MC,\bar q)$ the tuple $\bar p=\mu_\Delta(\bar{\bar p})$ is in $\varphi(\MB)$. So the $L(\MA)$-structure $\Gamma\circ\Delta(\MC,\bar r)=\Gamma\circ\Delta(\MC,(\bar {\bar p},\bar q))$ is well-defined and isomorphic to $\MA$, again due to item~\ref{le:int2}.  
\end{remark}

\subsection{Bi-interpretations and elementary equivalence}
\label{subsec:2.3}

In this section we discuss a very strong version of mutual interpretability of two structures, so-called {\em bi-interpretability}. The following definition uses notation from Lemma~\ref{le:int-transitivity}.

\begin{definition}\label{def:bi}
Algebraic structures $\MA$ and $\MB$ are called {\em strongly bi-interpretable} (with parameters) in each other, if 
there exists an interpretation $(\Gamma,\bar p,\mu_\Gamma)$ of $\MA$ into $\MB$ and an interpretation $(\Delta,\bar q,\mu_\Delta)$ of $\MB$ into $\MA$ (so the algebraic structures $\Gamma\circ\Delta(\MA,(\bar{\bar p},\bar q))$ and  $\Delta\circ\Gamma(\MB,(\bar{\bar q},\bar p))$  are uniquely defined 
 and $\Gamma\circ\Delta(\MA,(\bar{\bar p},\bar q))$ is isomorphic to $\MA$, while $\Delta\circ\Gamma(\MB,(\bar{\bar q},\bar p))$ is isomorphic to $\MB$), such that the coordinate maps $\mu_\Gamma\circ\mu_\Delta\colon A_{\Gamma \circ \Delta} \to A$ and $\mu_\Delta\circ\mu_\Gamma\colon B_{\Delta \circ \Gamma} \to B$ are definable in $\MA$ and $\MB$ respectively.
\end{definition}

Note that there is another slightly different notion of {\em bi-interpretation}, which we sometimes call a {\em weak bi-interpretation} for contrast, where in the above definition the condition of definability of maps $\mu_\Gamma\circ\mu_\Delta$ and $\mu_\Delta\circ\mu_\Gamma$ is replaced by a weaker one that requires definability of some coordinate maps $A_{\Gamma \circ \Delta} \to A$ and $B_{\Delta \circ \Gamma} \to B$. 

Often, the authors did not even mention the difference, implicitly assuming either one or another. To be precise, we endorse these two notions explicitly. Observe that the bi-interpretation defined in the books~\cite{Hodges} and~\cite{KhMS} is weak, but in the paper~\cite{AKNS} it is strong. There are many interesting applications of strong bi-interpretations which we cannot derive from the weak ones. However, it is important to mention that right now we do not have examples of weak interpretations of algebraic structures that are not strong.

\begin{definition}\label{def:reg}
Algebraic structures $\MA$ and $\MB$ are called {\em regularly bi-interpretable}, if  
\begin{enumerate}[label=\arabic*)]
\item there exist a regular interpretation $(\Gamma,\varphi)$ of $\MA$ in $\MB$ and a regular interpretation $(\Delta,\psi)$ of $\MB$ in $\MA$;
\item there exists formula $\theta_\MA(\bar u, x, \bar r)$ in $L(\MA)$, where $|\bar u|={\dim\Gamma\cdot\dim\Delta}$, $|\bar r|={\dim_{par}\Gamma\circ\Delta}$, such that for any tuple $\bar r_0\in \varphi_\Delta\wedge \psi(\MA)$ the formula $\theta_\MA(\bar u, x, \bar r_0)$ defines some coordinate map $U_{\Gamma\circ\Delta}(\MA,\bar r_0)\to A$;
\item there exists formula $\theta_\MB(\bar u, x, \bar t)$ in $L(\MB)$, where $|\bar u|={\dim\Gamma\cdot\dim\Delta}$, $|\bar t|={\dim_{par}\Delta\circ\Gamma}$, such that for any tuple $\bar t_0\in \psi_\Gamma\wedge \varphi(\MB)$ the formula $\theta_\MB(\bar u, x, \bar t_0)$ defines some coordinate map $U_{\Delta\circ\Gamma}(\MB,\bar t_0)\to B$. 
\end{enumerate}
\end{definition}

Regular bi-interpretability allows us to describe non-standard models of one algebraic structure by means of another one (see Theorem~\ref{th:equiv1} below).

\begin{definition}\label{def:st_reg}
We say, that $\MA$ and $\MB$ are {\em strongly regularly bi-interpretable}, if they are regularly bi-interpretable, i.\,e., 1)--3) hold, and additionally 
\begin{enumerate}
\item[4)] for any 
pair of parameters $(\bar p, \bar q)$, $\bar p\in \varphi(\MB)$, $\bar q\in \psi(\MA)$, there exists a pair of coordinate maps $(\mu_\Gamma,\mu_\Delta)$ for interpretations $(\Gamma,\bar p)$ and $(\Delta,\bar q)$, such that for any $\bar r_0=(\bar{\bar p},\bar q)$, $\bar{\bar p}\in \mu^{-1}_\Delta(\bar p)$, and $\bar t_0=(\bar{\bar q},\bar p)$, $\bar{\bar q}\in \mu^{-1}_\Gamma(\bar q)$, the coordinate maps ${\mu_\Gamma\circ\mu_\Delta\colon} U_{\Gamma\circ\Delta}(\MA,\bar r_0)\to A$ and ${\mu_\Delta\circ\mu_\Gamma\colon} U_{\Delta\circ\Gamma}(\MB,\bar t_0)\to B$ are defined in $\MA$ and $\MB$ correspondingly by the formulas $\theta_\MA(\bar u, x, \bar r_0)$ and $\theta_\MB(\bar u, x, \bar t_0)$.
\end{enumerate}
\end{definition}

\begin{remark}
Note that item~4) in Definition~\ref{def:st_reg} does not imply items~2), 3) in Definition~\ref{def:reg}. The statement in item~2) must be true for all tuples of parameters from the set $X=\varphi_\Delta\wedge\psi(\MA)$. From Remark~\ref{remark2} we know that  
$$
X=\{(\bar {\bar p},\bar q)\:\mid\: \bar q\in\psi(\MA), \:\bar p \in \varphi(\MB), \: \bar{\bar p}\in \mu_\Delta^{-1}(\bar p)\},
$$
where $\mu_\Delta$ runs all coordinate maps of the interpretation $\MB\simeq\Delta(\MA,\bar q)$. Denote by $\mu_\Delta^{\bar p,\bar q}$ the coordinate map of the interpretation $\MB\simeq\Delta(\MA,\bar q)$ from item~4), that corresponds to the pair of parameters $(\bar p,\bar q)$. Thus the statement in item~4) must be true for all tuples of parameters from the set 
$$
Y=\{(\bar {\bar p},\bar q)\:\mid\: \bar q\in\psi(\MA), \:\bar p \in \varphi(\MB),\: \bar{\bar p}\in (\mu_\Delta^{\bar p,\bar q})^{-1}(\bar p)\}.
$$
So $Y\subseteq X$, but it may be $Y\ne X$. Such example we will see in the Theorem~\ref{prop2} below.
\end{remark}

If algebraic structures $\MA$ and $\MB$ are strongly regularly bi-interpretable then we may use all conclusions from both regular bi-interpretability and strong bi-interpretation with parameters.

Recall that the \emph{first-order classification problem} for a structure $\MA$ asks one to describe "algebraically" all structures $\nsA$ such that $\MA \equiv \nsA$. 
In other words, the first-order classification problem for $\MA$ requires algebraically to describe all models of the complete first-order  theory $\Th(\MA)$ of $\MA$. Here, by "algebraically" we mean that the description has to reveal the algebraic structure of every model of $\Th(\MA)$.

The following result, in the case of regular bi-interpretation $\MA \simeq \Gamma(\MB,\phi)$, describes the algebraic structure of the models of $\Th(\MA)$ via the interpretation $\Gamma$ and the models of $\Th(\MB)$. The efficacy of this description depends on the understanding of the models $\nsB$ of $\Th(\MB)$. If  $\Gamma(\nsB,\phi)$ is well-defined, that is, the structures $\Gamma(\nsB,\bar p)$, $\bar p \in \phi(\nsB)$, are well-defined and pairwise isomorphic to each other, then we can choose an arbitrary representative $\bar p_0 \in \phi(\nsB)$ and view $\Gamma(\nsB,\phi)$ as $\Gamma(\nsB,\bar p_0)$. 

\begin{theorem}[\cite{Th_int2}]\label{th:equiv1}
Let $\MA$ and $\MB$ be regularly bi-interpretable in each other, so $\MA \simeq \Gamma(\MB,\phi)$ and $\MB\simeq\Delta(\MA,\psi)$. Then
\begin{enumerate}[label=(\arabic*)]
\item For any $\nsB\equiv\MB$ the algebraic structure $\Gamma(\nsB,\phi)$ is well-defined and $\MA\equiv \Gamma(\nsB,\phi)$;
\item Every $L(\MA)$-structure $\nsA$ which is elementarily equivalent to $\MA$ is isomorphic to  $\Gamma(\nsB,\phi)$ for some $\nsB\equiv\MB$;
\item For any $\MB_1\equiv\MB\equiv\MB_2$ one has 
$$
\Gamma(\MB_1,\phi) \simeq \Gamma(\MB_2,\phi) \iff \MB_1 \simeq \MB_2.
$$
\end{enumerate}
\end{theorem}


Following the practice in non-standard arithmetic and non-standard analysis, we term the structures of the form $\Gamma(\nsB,\bar p)$ for some $\nsB\equiv\MB$ and $\bar p \in \psi(\nsB)$ above as \emph{non-standard models} of $\MA$ with respect to interpretation $\Gamma$. In fact, the result below shows that quite often the non-standard models in $\nsB$ do not depend on the choice of interpretation $\Gamma$, i.e., there is only one up to isomorphism a non-standard model of $\MA$ in every $\nsB\equiv\MB$.

\begin{theorem} [Uniqueness of non-standard models, \cite{Th_int2}]\label{th:equiv2}
Let a finitely generated structure $\MA$ in a finite signature be regularly interpretable in $\Z$ in two ways, as $\Gamma_1(\Z,\phi_1)$ and $\Gamma_2(\Z,\phi_2)$. Then there exists a formula $\theta(\bar x_1,\bar x_2, \bar y_1,\bar y_2)$, $|\bar x_i|=\dim\Gamma_i$, $|\bar y_i|=\dim_{par}\Gamma_i$, such that for any $\bar p_i\in\phi_i(\Z)$ the formula $\theta(\bar x_1,\bar x_2, \bar p_1,\bar p_2)$ defines an isomorphism $\Gamma_1(\Z,\bar p_1)\to\Gamma_2(\Z,\bar p_2)$. Moreover, if $\nsZ \equiv \Z$, then algebraic structures $\Gamma_1(\nsZ,\phi_1), \Gamma_2(\nsZ,\phi_2)$ are well-defined and $\theta(\bar x_1,\bar x_2, \bar p_1,\bar p_2)$ defines an isomorphism $\Gamma_1(\nsZ,\bar p_1)\to\Gamma_2(\nsZ,\bar p_2)$ for any $\bar p_1\in\phi_1(\nsZ)$, $\bar p_2\in\phi_2(\nsZ)$.
\end{theorem}

\begin{corollary} \label{co:unique-non-stand}
    Let $\MA$ be a finitely generated structure in a finite signature regularly bi-interpretable with $\Z$. Then for every $\nsZ \equiv \Z$ there is a unique up to isomorphism non-standard model $\MA(\nsZ)$ of $\MA$ and for any structure $\nsB$ one has $\nsB \equiv \MA$ if and only if $\nsB \simeq \MA(\nsZ)$ for a suitable $\nsZ \equiv \Z$.
\end{corollary}

If $\MA$ is as in Corollary~\ref{co:unique-non-stand} then the structure $\MA(\nsZ)$ is called the  {\em non-standard model of $\MA$ with respect to $\nsZ$}. If $\MA \simeq \Gamma(\Z)$ is an interpretation of $\MA$ in $\Z$ then $\MA(\nsZ) \simeq \Gamma(\nsZ)$ and the algebraic structure of $\MA(\nsZ)$ is revealed via the interpretation $\Gamma$.
 
\section{Some facts from commutative algebra}
\label{se:3}

In this section, we discuss some results from commutative algebra that are certainly known in folklore. Since we were unable to find direct references, we provide short proofs.

Let $R$ be a commutative associative  unitary ring.  For $r, s \in R$ we write $r\mid s$ if $r$ divides $s$. By $R[x,x^{-1}]$ we denote the ring of Laurent polynomials in variable $x$ over $R$.  If $f\in R[x,x^{-1}]$ then
$\deg_{\min}(f)$ and $\deg_{\max}(f)$ denote the minimum and maximum powers of $x$ that occur in $f$ with non-zero coefficients. 

\begin{fact}\label{st1}
Let $R$ be a commutative associative unitary ring and $m, n\in \Z$, $n\ne 0$. Then
$$
n\mid m \mbox{ in the ring } \: \Z \iff (x^n-1)\mid (x^m-1) \mbox{ in the ring } R[x,x^{-1}].
$$
\end{fact}

\begin{proof}
For $m=0$ the fact is trivial, so we assume that $m\ne 0$. Suppose first that $m=n\cdot l$, for some $l\in \Z$. We need to show that there exists a polynomial $f\in R[x,x^{-1}]$ with 
$$
x^m-1=(x^n-1)\cdot f.
$$
Consider four cases:
\begin{itemize}
\item if $m,n>0$, then $l>0$ and
$$
x^m-1=(x^n)^l-1=(x^n-1)\cdot f_l(x^n),
$$
where $f_l(y)=1+y+\ldots+y^{(l-1)}$, $l>0$;
\item if $m, n<0$, then $l>0$, $(-m)=(-n)\cdot l$, and 
\begin{multline*}
x^m-1=-x^m(x^{-m}-1)=-x^m(x^{-n}-1)\cdot f_l(x^{-n})=x^{n+n(l-1)}(1-x^{-n})\cdot f_l(x^{-n})=\\=x^n(1-x^{-n})\cdot x^{n(l-1)}\cdot f_l(x^{-n})=(x^n-1)\cdot f_l(x^n);
\end{multline*}
\item if $m>0$, $n<0$, then $l<0$, $m=(-n)\cdot (-l)$, and
$$
x^m-1=(x^{-n}-1)\cdot f_{-l}(x^{-n})=- x^nx^{-n}(1-x^{-n})\cdot f_{-l}(x^{-n})=(x^n-1)\cdot g_l(x^n),
$$
where $g_l(y)=-y^{-1}-y^{-2}-\ldots-y^l$, $l<0$;
\item if $m<0, n>0$, then $l<0$, $(-m)=n\cdot (-l)$, and
$$
x^m-1=-x^m(x^{-m}-1)=-x^{nl}(x^n-1)\cdot f_{-l}(x^n)=(x^n-1)\cdot g_l(x^n).
$$
\end{itemize}

Assume now that $n\nmid m$, so $m=n\cdot l+r$, where $l,r\in \Z$ and $0<\vert r\vert<\vert n\vert$. This gives
$$
x^m-1=x^{nl+r}-1=x^{nl+r}-x^r+x^r-1=x^r(x^{nl}-1)+(x^r-1)\equiv x^r-1 (\mathrm{mod } (x^n-1)).
$$
It suffices to show that $(x^n-1)\nmid(x^r-1)$.  Suppose to the contrary that $(x^r-1)=(x^n-1)\cdot f$ for some polynomial $f\in R[x,x^{-1}]$.  Since $\deg_{\min}(f)=\deg_{\min}(x^r-1)-\deg_{\min} (x^n-1)$ and $\deg_{\max}(f)=\deg_{\max}(x^r-1)-\deg_{\max}(x^n-1)$, we obtain
\begin{itemize}
\item if $r, n>0$, then $\deg_{\min}(f)=0$ and $\deg_{\max}(f)=r-n$;
\item if $r, n<0$, then $\deg_{\min}(f)=r-n$ and $\deg_{\max}(f)=0$;
\item if $r>0, n<0$, then $\deg_{\min}(f)=-n$ and $\deg_{\max}(f)=r$;
\item if $r<0, n>0$, then $\deg_{\min}(f)=r$ and $\deg_{\max}(f)=-n$.
\end{itemize}
In each case we get a contradiction with  $\deg_{\min}(f)\leqslant\deg_{\max}(f)$ and $0<\vert r\vert<\vert n\vert$.
\end{proof}

Recall that for any $k\in \N$, $k>1$, 
$\Z[1/k] =\{zk^i \mid z,i\in\Z\}$
 is a subring of the ring of rationals $\Q$, containing  $\Z$.

\begin{fact}\label{cor3}
Let $k\in \N$, $k>1$, and $m, n\in \Z$, $n\ne 0$. Then
$$
n\mid m \mbox{ in the ring } \: \Z \iff (k^n-1)\mid (k^m-1) \mbox{ in the ring } \: \Z[1/k].
$$
\end{fact}

\begin{proof}
If $n\mid m$, then by Fact 1 there exists a polynomial $f\in \Z[x,x^{-1}]$ with $x^m-1=(x^n-1)\cdot f$. Specifying $x$ into $k$ one gets $k^m-1=(k^n-1)\cdot f(k)$, so $(k^n-1)\mid (k^m-1)$ in the ring $\Z[1/k]$. 

On the other hand, let $m=n\cdot l+r$, $l,r\in \Z$, and $0<\vert r\vert<\vert n\vert$, but $(k^r-1)=(k^n-1)\cdot z$ for an element $z\in\Z[1/k]$. We can write $z=z_0\cdot k^{-i}$, where $z_0, i\in \Z$, $i\geqslant 0$. 

If $i=0$ then $(k^r-1)=(k^n-1)\cdot z_0$ and
\begin{itemize}
\item if $r,n>0$, the this  contradicts to $\vert k^r-1\vert<\vert (k^n-1)\cdot z_0 \vert$;
\item if $r, n<0$, then $k^{\vert n\vert -\vert r \vert}\cdot(1-k^{\vert r\vert})=(1-k^{\vert n\vert})\cdot z_0$ and since $(k,(k^{\vert n\vert}-1))=1$, we obtain that $k^{\vert n\vert -\vert r \vert}\mid z_0$ in $\Z$ and find a contradiction as above;
\item the cases  $r>0, n<0$ and $r<0, n>0$ are similar to the one above.
\end{itemize}

Now, let $i>0$ and $(z_0,k)=1$. Thus, $k^i\cdot(k^r-1)=(k^n-1)\cdot z_0$ and $0<\vert r\vert<\vert n\vert$. If $n<0$ then $k^{\vert n\vert}\cdot k^i\cdot(k^r-1)=(1-k^{\vert n\vert})\cdot z_0$ and hence $k\mid z_0$. Now, assume $n>0$. If $r<0$, then $k^i \cdot (1-k^{\vert r\vert})=(k^n-1)\cdot z_0\cdot k^{\vert r\vert}$. If $\vert r\vert<i$ then $k\mid z_0$, hence we may assume that $\vert r\vert\geqslant i$.  If $r <0$ then $k^i\cdot(k^r-1)=(k^n-1)\cdot z_0$ does not hold, since $z_0$ is an integer.  Therefore,  $r, n, i>0$ and $z_0=z_0\cdot k^n-(k^r-1)\cdot k^i$, so we obtain a contradiction with $(z_0,k)=1$ again. 
\end{proof}

\begin{fact}\label{cor1}
Let $n,l\in \Z$, then there exists a polynomial $g\in \Z[x,x^{-1}]$ such that 
$$
x^{nl}-1=(x^n-1)\cdot (l+(x^n-1)\cdot g).
$$
\end{fact}

\begin{proof}
If $n=0$ or $l=0$ then the statement is obvious, so we assume that $n,l\ne 0$.

If $l>0$ then $x^{nl}-1=(x^n-1)\cdot f_l(x^n)$, where $f_l(y)=1+(y-1+1)+(y-1+1)^2+\ldots+(y-1+1)^{l-1}=l+(y-1)\cdot h_l$ for some polynomial $h_l\in \Z[y]$. In this case, the desired polynomial $g$ equals $h_l(x^n)$.

If $l<0$ then $x^{nl}-1=(x^n-1)\cdot g_l(x^n)$, where 
$$
g_l(y)=-y^{-1}-y^{-2}-\ldots-y^l=-y^l\cdot f_{-l}(y)=-y^l\cdot(-l+(y-1)\cdot h_{-l})=y^l(l-(y-1)\cdot h_{-l}).
$$
Thus we have
\begin{multline*}
x^{nl}-1=(x^n-1)\cdot x^{nl}(l-(x^n-1)\cdot h_{-l}(x^n))=(x^n-1)\cdot(x^{nl}-1+1)\cdot (l-(x^n-1)\cdot h_{-l}(x^n))=\\
=(x^n-1)\cdot (l-(x^n-1)\cdot h_{-l}(x^n)+(l-(x^n-1)\cdot h_{-l}(x^n))\cdot(x^{nl}-1)).
\end{multline*}
And since $x^{nl}-1=(x^n-1)\cdot g_l(x^n)$ it suffices to  put
$$
g=(l-(x^n-1)\cdot h_{-l}(x^n))\cdot g_l(x^n) - h_{-l}(x^n).
$$
\end{proof}

\begin{corollary}\label{cor2}
Let $z\in \Z$ and $k\in \N$, $k>1$. Then for every $n\in \Z\setminus\{0\}$ one has  
$$
\frac{k^{nz}-1}{k^n-1} \, \equiv\,  z \,(\mbox{\rm mod} \,(k^n-1)) 
$$
 in the ring  $\Z[1/k]$.
\end{corollary}

\begin{fact}\label{fact1}
Let $k\in \N$, $k>1$, and $y\in \Z[1/k]$ such that for every $n\in \Z^+$ one has  
$$
y \, \equiv\,  0 \,(\mbox{\rm mod} \,(k^n-1)) 
$$
in the ring  $\Z[1/k]$. Then $y=0$.
\end{fact}

\begin{proof}
Assume that $y=y_0 \cdot k^{-i}$, $y_0\in \Z$, $i\in\Z$, $i\geqslant 0$. For every $n\in \Z^+$ there exists $g_n\in \Z[1/k]$ such that $y_0=(k^n-1)\cdot g_n\cdot k^i$. As $y_0\in \Z$ and $(k^n-1,k)=1$, it should be $g_n\cdot k^i\in \Z$. Thus $(k^n-1)\mid y_0$ in $\Z$ for every $n\in \Z^+$. But $k^{2n}-1=(k^n-1)\cdot (k^n+1)$ and $\mathrm{GCD}((k^n-1),(k^n+1))\in\{1,2\}$, therefore the set of distinct prime divisors for numbers $\{k^n-1, n\in \Z^+\}$ is infinite. It implies that $y_0=0$, hence $y=0$, as required.
\end{proof}

\begin{corollary}\label{cor4}
Let $k\in \N$, $k>1$, and $z\in \Z$, $l,t\in\Z[1/k]$. Then the following conditions are equivalent:
\begin{enumerate}
\item $l\cdot z=t$;
\item $l\cdot(k^{nz}-1)  \equiv  t\cdot (k^n-1) \:(\mathrm{mod } \,(k^n-1)^2)$ in the ring $\Z[1/k]$ for every $n\in \Z$.
\end{enumerate}
\end{corollary}

\begin{proof}
For $(1)\Longrightarrow(2)$ we use Corollary~\ref{cor2}:
\begin{multline*}
(1) \;\Longrightarrow\;  l\cdot z \equiv t \: (\mbox{\rm mod} \,(k^n-1))\; \forall\: n\in \Z \;\Longrightarrow \\ \Longrightarrow\; l\cdot(k^{nz}-1)/(k^n-1)\equiv t \: (\mbox{\rm mod} \,(k^n-1)) \;\forall\: n\in \Z\setminus\{0\} \;\Longrightarrow \; (2).
\end{multline*}
And for $(2)\Longrightarrow(1)$ we refer to Fact~\ref{fact1} and Corollary~\ref{cor2}:
\begin{multline*}
(2) \;\Longrightarrow l\cdot(k^{nz}-1)/(k^n-1)-t\equiv 0 \: (\mbox{\rm mod} \,(k^n-1)) \;\forall\: n\in \Z^+ \;\Longrightarrow\\ \;\Longrightarrow \; l\cdot z -t \equiv 0 \: (\mbox{\rm mod} \,(k^n-1))\; \forall\: n\in \Z^+ \;\Longrightarrow\; (1).
\end{multline*}
\end{proof}

As usual, by $U(R)$ we denote the multiplicative group of units (invertible elements) in a ring $R$. Note that 
$U(\Z[1/k])$ is generated (as a multiplicative group) by all divisors of $k$, i.\,e., 
$$
U(\Z[1/k]) = \langle d \mid d \ \mathrm{divides} \ k\rangle.
$$

\begin{fact}\label{fact2}
Let $k\in \N$, $k>1$, and $y\in \Z[1/k]$. Then one has
$$
y\in U(\Z[1/k]) \;\Longleftrightarrow 
\;\; \forall \: t\in \Z[1/k] \;\; \exists \: z,i\in \Z,  \;\;  \quad y\cdot z= t\cdot k^i.
$$
\end{fact}

\begin{proof}
If $y\in U(\Z[1/k])$ then there exists an element $1/y\in \Z[1/k]$. For any $t\in Z[1/k]$ we can f
ind $i\in \Z$ such that $z=t\cdot 1/y \cdot k^i\in \Z$ and the result follows.

Suppose now that $y$ satisfies the converse condition. Take $t=1$. There exist $z,i\in \Z$  such that $y\cdot z=k^i$, hence $1/y=z\cdot k^{-i}$.
\end{proof}

\section{Regular bi-interpretation of $BS(1,k)$ and $\Z$}
\label{subsec:4.1}

We start with reminding some known facts about the metabelian Baumslag\,--\,Solitar groups $BS(1,k)$. As was mentioned above the group $BS(1,k)$ is  defined by one-relator presentation:
$$
BS(1,k)=\langle a,b \mid b^{-1}ab=a^k \rangle,
$$
where $k \in \Z^+$. If $k = 1$, then $BS(1,1)$ is free abelian of rank $2$. To exclude this trivial case we always assume that $k \geq 2$. Clearly, $a^zb^i=b^ia^{zk^i}$ (or $b^{-i}a^z=a^{zk^i}b^{-i}$) for all $z,i\in\Z$, $i\geqslant 0$, thus any element $g\in BS(1,k)$ has a form $g=b^i\,a^z\,b^{-i}\,b^m$, $z,i, m\in \Z$, $i\geqslant 0$. 

\subsection{$BS(1,k)$ as a semi-direct product}
It is known  that the group $BS(1,k)$ is isomorphic to the semi-direct product $\Z[1/k]\rtimes \Z$ of abelian groups, where 
$$
\Z[1/k] =\{zk^i\mid  z,i\in\Z\}
$$ 
is a subgroup (in fact, a subring) of the additive group $\langle \Q; +,-,0 \rangle$ of rationals (the ring $\Q$), and the action $\varphi\colon \Z \to \Aut(\Z[1/k])$ of $\Z$ on $\Z[1/k]$ is given by 
$$
\varphi(m)=\varphi_m,\quad \varphi_m(x) = x\cdot k^{-m}, \quad x\in\Z[1/k],\: m\in \Z.
$$

Recall that elements in $\Z[1/k]\rtimes \Z$ are pairs $(zk^{i},m)$, where $zk^{i}\in\Z[1/k]$, $m\in \Z$; the identity element is $(0,0)$, and the product in $\Z[1/k]\rtimes \Z$ is defined as
\begin{equation}  \label{eq:product}
(z_1k^{i_1}, m_1)(z_2k^{i_2}, m_2)=(z_1k^{i_1}+z_2k^{i_2-m_1}, m_1+m_2).
\end{equation}
The inverse of an element $(zk^{i},m)$ is defined by 
\begin{equation} \label{eq:inverse}
(zk^{i},m)^{-1} = (-zk^{i+m}, -m). 
\end{equation}
It is follows from~\eqref{eq:product} and~\eqref{eq:inverse} that 
for $x_1=(y_1,m_1)$, $x_2=(y_2,m_2)$, $y_i\in \Z[1/k]$, $m_i\in \Z$, one has
\begin{gather*}
x_2^{-1}x_1x_2=(y_2,m_2)^{-1}(y_1,m_1)(y_2,m_2)=((y_1-y_2)k^{m_2}+y_2k^{m_2-m_1},m_1),\\ 
[x_1,x_2]=[(y_1,m_1),(y_2,m_2)]=(-y_1k^{m_1}+y_2k^{m_2}+(y_1-y_2)k^{m_1+m_2},0),
\end{gather*}
in particular, if $m_1=0$, then 
\begin{equation}\label{eq:semi}
x_2^{-n}x_1x_2^n=(y_2,m_2)^{-n}(y_1,0)(y_2,m_2)^n=(0,m_2)^{-n} x_1 (0,m_2)^n, \quad n\in \Z.
\end{equation}

For elements $a=(1,0)$ and $b=(0,1)$ one has $b^{-i} a^z\, b^{i}\, b^m=(zk^{i},m)$, $z,i,m\in \Z$. Therefore, $a$ and $b$ generate the group $\Z[1/k]\rtimes \Z$ and $b^{-1}a\,b=a^k$. Thus the map $a \to (1,0)$, $b \to (0,1)$ gives rise to an isomorphism  $\lambda\colon BS(1,k)\to \Z[1/k]\rtimes \Z$. Note also that $\Z[1/k] \cong \ncl(a)=C_{BS(1,k)}(a)$ and $\Z \cong \langle b \rangle=C_{BS(1,k)}(b)$. 

Furthermore, the subgroup $\ncl(a)$ has a structure of $\Z[1/k]$-module and it is torsion-free. It is convenient to use powers when writing multiplication of elements $x=(y,0)$ from $\ncl(a)$ to coefficients from $\Z[1/k]$: 
$$
x^{z k^{i}}=b^{-i} x^z b^{i}=(y,0)^{z k^{i}}=(y\cdot zk^i,0), \quad y\in\Z[1/k], \; z,i\in \Z.
$$
It gives an easy form for elements from $\ncl(a)$, namely, $\ncl(a)=\{a^y\,\vert\, y\in\Z[1/k]\}$. Therefore, every element $g\in BS(1,k)$ has a form $g=a^y\,b^m$, where $y\in \Z[1/k]$, $m\in \Z$ are uniquely defined:  $g=(y,m)=a^{y}b^m$. 

If $m\ne 0$ then by induction on $\vert n\vert$ it easy to see, that 
\begin{equation}\label{eq:g^n}
g^n=(y,m)^n=(y\frac{k^{-mn}-1}{k^{-m}-1},mn), \quad n\in \Z.
\end{equation}

\subsection{Absolute interpretation of $BS(1,k)$ in $\Z$}
\label{subsec:4.2}

Let us define an interpretation $\Delta$ of the algebraic structure $\langle BS(1,k), \,\cdot \, , \, ^{-1}, e \rangle$ in the algebraic structure $\langle \Z, +, \, \cdot\, , 0, 1 \rangle$ with the base $\Z^3$ and the coordinate map
\begin{equation}\label{eq:mu_Delta}
\mu_\Delta\colon (z,i,m)\in\Z^3 \;\longrightarrow \; (zk^{i},m)\in BS(1,k).
\end{equation}
An equivalence on $\Z^3$ is defined by the rule
\begin{equation}\label{eq:sim_Delta}
(z_1,i_1,m_1)\sim_\Delta(z_2,i_2,m_2) \iff z_1k^{i_1}=z_2k^{i_2} \; \wedge \; m_1=m_2.
\end{equation}
The rules~~\eqref{eq:product}, \eqref{eq:inverse}, \eqref{eq:sim_Delta}, that determine the relation $\sim_\Delta$  and group operations  on triples, allow us to assert that they are definable in $\Z$. And also $\Z^3/{\sim_\Delta}$ is isomorphic to $BS(1,k)$. Thus $BS(1,k)$ is absolutely interpretable in $\Z$.

To construct a converse interpretation $\Gamma$ of $\Z$ in $BS(1,k)$ we need some information about systems of generators in the group $BS(1,k)$.

\subsection{Another generators of $BS(1,k)$}

Along with the standard generators $a,b$ of the group $BS(1,k)$, it is convenient for us to consider other pairs of generators. We will use the following notations:
\begin{gather*}
A=\ncl(a)=\{a^y\in A \mid y\in \Z[1/k]\},\\  
A_1=\{a^y\in A \mid y\in U(\Z[1/k])\},\\
Ab=\{a^y b \mid y \in \Z[1/k]\}.
\end{gather*}

Note that an element $a_1\in BS(1,k)$ belongs to $A_1$ if and only if $\ncl(a_1)=A$. Also an element $b_1\in BS(1,k)$ belongs to $Ab$ if and only if $b_1=(y,1)$, $y\in\Z[1/k]$. Therefore, since~\eqref{eq:semi}, all elements form $Ab$ act on $A$ by the same way as $b$ does:
\begin{equation}\label{lemma4}
b_1^{-n}\,u\, b_1^n=b^{-n}\, u \, b^n=u^{k^n}, \quad b_1\in Ab, \: u\in A, \: n\in \Z.
\end{equation}

\begin{lemma}\label{lemma5}
Let $a_1\in A_1$ and $b_1 \in Ab$.  Then the map $a\to a_1$, $b\to b_1$ defines an automorphism $\lambda$ of the group $BS(1,k)$.
\end{lemma}

\begin{proof}
By~\eqref{lemma4}, one has $b_1^{-1}a_1 b_1=a_1^k$, therefore, there exists an endomorphism $\lambda$ of $BS(1,k)$, such that $\lambda(a)=a_1$, $\lambda(b)=b_1$, and $\lambda(a^y\,b^m)=a_1^y\,b_1^m$ for all $y\in \Z[1/k]$, $m\in \Z$. Since $a_1^y\,b_1^m=(*,m)$, then $a_1^y\,b_1^m=e$ implies $m=0$ and after that $y=0$, so $\lambda$ is a monomorphism.

Suppose that $a_1=a^z$ with $z\in U(\Z[1/k])$, so  $1/z\in\Z[1/k]$. Then for an arbitrary $a^y\in A$, $y\in \Z[1/k]$, one has 
$a^y= a^{y\cdot z\cdot 1/z}=a_1^{y\cdot 1/z}=\lambda(a^{y\cdot 1/z})$, therefore, $\ncl(a)\subset\lambda(BS(1,k))$. Since $b_1=a_0 b$, $a_0\in A$, then we get $b=a_0^{-1}b_1\in \lambda(BS(1,k))$.
So $\lambda$ is an epimorphism and, therefore, it is an automorphism. 
\end{proof}

\begin{corollary}\label{cor6}
If $a_1\in A_1$ and $b_1\in Ab$, then every element $g\in BS(1,k)$ has a form $g=a_1^y\,b_1^m$, where $y\in Z[1/k]$, $m\in \Z$ are uniquely defined.
\end{corollary}

\begin{corollary}\label{lemma3}
For an element $b_1\in Ab$ one has $C_{BS(1,k)}(b_1)=\langle b_1 \rangle$.
\end{corollary}

It is important for our purposes to show that the sets $A, Ab, A_1$ are $0$-definable in $BS(1,k)$. We prove it in Lemmas~\ref{lemma1}, \ref{lemma2}, \ref{lemma8} below. By the way, it follows that the group $BS(1,k)$ has no automorphisms other than those described in Lemma~\ref{lemma5}.

\begin{lemma}\label{lemma1}
The formula
$$
\alpha(x)\;=\;\forall\: y  \;([y^{-1}xy,x]=e)
$$
defines $A$ in $BS(1,k)$.
\end{lemma}

\begin{proof}
If $x\in \ncl(a)$, then $\alpha(x)$ is obviously true in $BS(1,k)$ on $x$. Conversely, let $x=(y,m)$, $y\in\Z[1/k]$, $m\in \Z$, then
$$
h=a^{-1}xa=(1,0)^{-1}(y,m)(1,0)=(y-1+k^{-m},m),
$$
$$
g=[h,x]=[(y-1+k^{-m},m),(y,m)]=(-(k^m-1)^2,0),
$$
so if $m\ne 0$, i.\,e., $x\notin \ncl(a)$, then $g\ne (0,0)$ and $\alpha(x)$ is false in $BS(1,k)$ on $x$.    
\end{proof}

\begin{lemma}\label{lemma2}
The formula
$$
\beta(y)\;=\; \forall \: x\;\: (\alpha(x)\;\longrightarrow \; y^{-1}xy=x^k)
$$
defines $Ab$ in $BS(1,k)$.
\end{lemma}

\begin{proof}
Take an element $y=(x,1)\in Ab$, $x\in\Z[1/k]$. For any element $h=(z,n)$, $z\in \Z[1/k]$, $n\in\Z$, one has
$$
g=y^{-1}hy=(x,1)^{-1}(z,n)(x,1)=((z-x)k+xk^{1-n},n),
$$
in particularly, if $(z,n)\in \ncl(a)$, i.\,e., $n=0$, then $g=(zk,0)=(z,0)^k$. Thus any element $(x,1)$ satisfies the formula $\beta$. Conversely, let $w=(x,m)$, $x\in\Z[1/k]$, $m\in \Z$, then  
$$
g=w^{-1}aw=(x,m)^{-1}(1,0)(x,m)=(k^m,0)=(1,0)^{k^m}, 
$$
so if $w\not\in Ab$, i.\,e., $m\ne 1$, then $g\ne (1,0)^k$, thus the element $w$ does not satisfy the formula $\beta$. 
\end{proof}

\subsection{Regular interpretation of $\Z$ in $BS(1,k)$}

Now let us define a regular injective interpretation $(\Gamma, \beta)$ of $\Z$ in $BS(1,k)$ with $\dim\Gamma=\dim_{par}\Gamma=1$. First we put
$$
U_\Gamma(x,b_1)\:=\: ([x,b_1]=e),
$$
and $\sim_\Gamma$ is the identity $=$ on the set $U_\Gamma(BS(1,k),b_1)$.
Take any parameter $b_1\in \beta(BS(1,k))$. By Corollary~\ref{lemma3} and Lemma~\ref{lemma2}, one has $U_\Gamma(BS(1,k),b_1)=\langle b_1\rangle$. We define the coordinate map in this way:
\begin{equation}\label{eq:mu_Gamma}
\mu_\Gamma\colon b_1^m \in \langle b_1\rangle \; \longrightarrow\; m\in \Z.
\end{equation}
It is obvious that addition of integers is definable on $\langle b_1\rangle$ with parameter $b_1$, also $0$ corresponds to $e$ and $1$ corresponds to $b_1$. It remains only to show that multiplication of integers is definable on $\langle b_1\rangle$, namely, that there exists a formula $\gamma(x,y,z,t)$ in the group language, such that 
\begin{equation}\label{eq:multiplication}
BS(1,k)\models \gamma(b_1^n,b_1^l,b_1^m,b_1) \iff n\cdot l=m, \quad n,l,m\in \Z.  
\end{equation}

\begin{lemma}\label{lemma6}
The multiplication of integers is definable on $\langle b_1\rangle$ with parameter $b_1$, provided $b_1\in Ab$.
\end{lemma}

\begin{proof}
To prove that multiplication is definable on $\langle b_1\rangle$ it is enough to show that divisibility is definable on $\langle b_1\rangle$ as well~\cite{Robinson}, because 
\begin{multline}\label{eq:div}
m=n\cdot l \iff (n+l)\cdot (n+l+1)=n\cdot(n+1)+l\cdot(l+1)+m+m,\\
l=n\cdot (n+1) \iff \\ \iff (n+n+1)\mid (l+l-n)\:\wedge \\ \wedge\:(\forall \:m\;(l\mid m\longleftrightarrow (n\mid m \:\wedge\: (n+1)\mid m))), \quad m,n,l\in\Z\setminus\{0,-1\}.
\end{multline}

According to Fact~\ref{cor3}, $n\mid m$ in $\Z$ if and only if $(k^n-1)\mid (k^m-1)$ in $\Z[1/k]$. The last means that there exits $t\in \Z[1/k]$ such that $k^m-1=t\cdot(k^n-1)$. It gives that
$c^{k^m-1}=(c^t)^{k^n-1}$ for any $c\in \ncl(a)$, i.\,e., $b^{-m}cb^{m}c^{-1}=b^{-n}c^tb^nc^{-t}$. Let us define the following formula: 
$$
\delta(x,y)\;= \;\forall \: c\; (\alpha(c)\;\longrightarrow\; \exists\: u\; (\alpha(u)\wedge [x,c]=[y,u])).
$$
Suppose that $x,y\in \langle b_1\rangle$, say $x=b_1^m, y=b_1^n$, and $n\ne 0$. If $n\mid m$, then for any $c\in \alpha(BS(1,k))=\ncl(a)$ there exists $u=c^t$, such that $b^{-m}c^{-1}b^{m}c=b^{-n}u^{-1}b^nu$, and by~\eqref{lemma4}, $b_1^{-m}c^{-1}b_1^{m}c=b_1^{-n}u^{-1}b_1^nu$, i.\,e., $[x,c]=[y,u]$. Therefore, $BS(1,k)\models \delta(x,y)$. Inversely, if $BS(1,k)\models \delta(x,y)$, then for $c=a^{-1}$ we obtain $a^{k^m-1}=u^{k^n-1}$ for some $u\in A$, i.\,e., there exists $t\in \Z[1/k]$, such that $u=a^t$, and, therefore, $k^m-1=t\cdot(k^n-1)$, thus $n\mid m$. So the formula $\delta(x,y)$ defines divisibility of integers on $\langle b_1\rangle$.

Now one can construct the formula $\gamma$~\eqref{eq:multiplication} for multiplication of integers on $\langle b_1\rangle$ using the formula $\delta$ and formulas~\eqref{eq:div}, and separately describing cases $n,l,m\in\{-1,0\}$. 
\end{proof}

\subsection{Strong regular bi-interpretation of $BS(1,k)$ and $\Z$}
\label{subsec:4.5}

\begin{theorem} \label{prop2}
If $k>1$, then the group $BS(1,k)$ is regularly strongly bi-interpretable with the ring $\Z$.
\end{theorem} 

Notice that we have defined below the regular injective interpretation $(\Gamma,\beta)$ of $\Z$ in $BS(1,k)$ and the absolute interpretation $\Delta$ of $BS(1,k)$ in $\Z$. And now we need to show that there exist formulas $\theta_\Z$ and $\theta_{BS(1,k)}$, that define coordinate maps for the compositions $\Gamma\circ\Delta$ and $\Delta\circ\Gamma$. We start here with several preliminary results.

\begin{lemma}\label{lemma7}
The formula
\begin{multline*}
\tau(x,y,h,b_1)\;= \;
\alpha(x)\,\wedge\,\alpha(y) \,\wedge\,([h,b_1]=e)\,\wedge\,\beta(b_1)\,\wedge\\
\wedge\,\forall\:v\;\forall \: w\;([v,b_1]=[w,b_1]=e \,\wedge\, \gamma(h,v,w,b_1)\;\longrightarrow \;\\\longrightarrow \;\exists\:u\;(\alpha(u)\,\wedge\,[v,y]=[w,x][v,[u,v]])\,)
\end{multline*}
is true in $BS(1,k)$ on elements $x,y,h,b_1$ if and only if $x,y\in A$ (i.\,e., $x=a^l$, $y=a^t$, $t,l\in \Z[1/k]$), $b_1 \in Ab$, $h\in\langle b_1\rangle$ (i.\,e., $h=b_1^z$, $z\in \Z$) and $l\cdot z=t$.
\end{lemma}

\begin{proof}
Let $z\in\Z$ and $t,l\in \Z[1/k]$. According to Corollary~\ref{cor4}, the identity $l\,\cdot\, z=t$ takes place if and only if for every $n\in \Z$ there exists $s_n\in\Z[1/k]$, such that 
\begin{equation}\label{eq:s_n}
 l\cdot(k^{nz}-1)=t\cdot(k^n-1)+s_n(k^n-1)^2.   
\end{equation}
For any $u\in A$ we may write $u^{k^n-1}=u^{k^n}u^{-1}=b^{-n}ub^nu^{-1}=[b^n,u^{-1}]$ and, by~\eqref{lemma4}, $[b^n,u^{-1}]=[b_1^n,u^{-1}]$. So let us continue: 
\begin{multline*}
\eqref{eq:s_n} \iff -t\cdot(k^n-1)=-l\cdot(k^{nz}-1)+s_n(k^n-1)^2 \iff \\ \iff (a^{-t})^{k^n-1}=(a^{-l})^{k^{nz}-1}(a^{s_n})^{(k^n-1)^2} \iff \\ \iff [b^n,a^t]=[b^{nz},a^l][b^n,[a^{-s_n},b^n]] \iff [b_1^n,a^t]=[b_1^{nz},a^l][b_1^n,[a^{-s_n},b_1^n]].
\end{multline*}

If $x=a^l$, $y=a^t$, $h=b_1^z$ and $l\cdot z=t$, then $BS(1,k)\models\tau(x,y,h,b_1)$. Here for any elements $v, w$ if $BS(1,k)\models([v,b_1]=[w,b_1]=e) \,\wedge\, \gamma(h,v,w,b_1)$, then there exists $n\in \Z$ such that $v=b_1^n$ and $w=b_1^{nz}$, thus one can take $u=a^{-s_n}$.

Inversely, if $x,y,h,b_1$ are elements from the group $BS(1,k)$ and  $BS(1,k)\models\gamma(x,y,h,b_1)$, then $b_1\in Ab$ and $x=a^l$, $y=a^t$, $h=b_1^z$ for some $z\in\Z$, $t,l\in \Z[1/k]$. For any $n\in \Z$ there exits $u\in A$ such that $[b_1^n,a^t]=[b_1^{nz},a^l][b_1^n,[u,b_1^n]]$. Let $u=a^f$, $f\in \Z[1/k]$, then for $s_n=-f$ one has~\eqref{eq:s_n}. Therefore, $l\cdot z=t$.
\end{proof}

\begin{corollary}\label{cor5}
If $a_1\in A$ and $b_1\in Ab$, then
$$
BS(1,k)\models\tau(a_1,u,h,b_1)\iff \exists\,z\in \Z \;(u=a_1^z \wedge h=b_1^z).
$$
\end{corollary}

\begin{proof}
Let $a_1=a^l$ and $l\in \Z[1/k]$.
If $u=a_1^z$ and $h=b_1^z$, $z\in \Z$, then $\tau(a_1,u,h,b_1)=\tau(a^l,a^{lz},b_1^z,b_1)$, therefore, $BS(1,k)\models\tau(a_1,u,h,b_1)$. Inversely, if $BS(1,k)\models\tau(a_1,u,h,b_1)$, then there exist $z\in\Z$ and $t\in\Z[1/k]$, such that $u=a^t$, $h=b_1^z$ and $l\cdot z=t$. Therefore, $u=a^{lz}=a_1^z$, as required.
\end{proof}

\begin{lemma}\label{lemma8}
The formula
$$
\pi(x)\;=\; \forall \: c \;\forall \: v\; (\alpha(c)\,\wedge\,\beta(v)\;\longrightarrow \;\exists \: w\;\exists\:h\; ([w,v]=e\,\wedge\,\tau(x,w^{-1}qw,h,v)\,)\,)
$$
defines the subset $A_1$ in $BS(1,k)$.
\end{lemma}

\begin{proof}
Indeed, Fact~\ref{fact2} says that $y\in U(\Z[1/k])$ if and only if for any $t\in \Z[1/k]$ there exist $z,i\in \Z$, such that $y\cdot z=t\cdot k^i$. By Lemma~\ref{lemma7}, for all  $t\in \Z[1/k]$, $i,z\in\Z$ one has 
$$
y\cdot z=t\cdot k^i \iff BS(1,k)\models\tau (a^y,a^{tk^i},b_1^z,b_1)\quad\forall (\exists)\:b_1\in Ab.
$$
Since $a^{tk^i}=b^{-i}\,a^t\,b^i=(b_1^i)^{-1}\,a^t\,b_1^i$, we obtain that $a^y\in A_1$ if and only if for any $a^t\in A$ and any $b_1\in Ab$ there exist $b_1^i,b_1^z\in \langle b_1\rangle$, such that $BS(1,k)\models\tau(a^y,\,(b_1^i)^{-1}\,a^t\,b_1^i, \,b_1^z, \,b_1)$. Therefore, $x\in A_1$ if and only if $BS(1,k)\models\pi(x)$.

\end{proof}

\begin{proof}[Proof of Theorem~\ref{prop2}] 
We have to show the existence of  formulas $\theta_\Z$ and $\theta_{BS(1,k)}$ that define coordinate maps
for the compositions $\Delta\circ\Gamma$ and $\Gamma\circ \Delta$. Let us remind that the interpretation $\Delta$ is absolute, while $\dim\Delta=3$, and $\Gamma$ has one parameter $b_1$, while $\dim\Gamma=1$. However, we may assume that $\Gamma$ has two parameters $a_1\in \pi(BS(1,k))$ and $b_1\in \beta(BS(1,k))$. We add to formulas from the code $\Gamma$ new free variable $\tilde a$ for parameter $a_1$, but in actually they do not depend on variable $\tilde a$. In this case $\Gamma$ remains regular with the formula $\phi(\tilde a,\tilde b)=\pi(\tilde a)\wedge\beta(\tilde b)$. 

So the compositions of the interpretations $(\Gamma,\phi)$ and $(\Delta,\emptyset)$ have forms $(\Gamma\circ\Delta,\phi_\Delta)$ and $(\Delta\circ\Gamma,\phi)$. Here formula $\phi_\Delta$ depends on six free variables ($z_a,i_a, m_a,z_b,i_b,m_b$). 
Formulas $\theta_\Z$ and $\theta_{BS(1,k)}$ that we are looking for have four free variables $(\tilde z,\tilde \imath,\tilde m, \,x)$, and additionally $\theta_{BS(1,k)}$ has two variables for parameters, while $\theta_\Z$ has six. Let us define the formulas $\theta_\Z$, $\theta_{BS(1,k)}$ and then check that they meet all the requirements:
\begin{gather*}
\theta_{BS(1,k)}(\tilde z,\tilde\imath,\tilde m,\:x, \:\tilde a,\tilde b)\;=\;([\tilde\imath,\tilde b]=[\tilde m,\tilde b]=e)\:\wedge\:\tau(\tilde a,\,\tilde\imath \, x \, \tilde m^{-1} \tilde\imath^{-1},\,\tilde z,\,\tilde b),
\\
\theta_{\Z}(\tilde z,\tilde \imath, \tilde m,\:x,\:z_a,i_a,m_a,\,z_b,i_b,m_b)\;=\;(x=\tilde m)\wedge (\tilde z k^{\tilde\imath}=z_bk^{i_b}\frac{k^{-\tilde m}-1}{k^{-1}-1}).
\end{gather*}
Remind the requirements for formulas $\theta_\Z$ and $\theta_{BS(1,k)}$ from Definition~\ref{def:reg}:
\begin{enumerate}
\item\label{req:1} for any parameters $\bar r\in \phi_\Delta(\Z)$ and $(a_1,b_1)\in\phi(BS(1,k))$ there exist coordinate maps $U_{\Gamma\circ\Delta}(\Z,\bar r)\to \Z$ and $U_{\Delta\circ\Gamma}(BS(1,k),a_1,b_1)\to BS(1,k)$ that are definable by formulas $\theta_\Z(\tilde z,\tilde \imath, \tilde m,\:x,\:\bar r)$ and $\theta_{BS(1,k)}(\tilde z,\tilde\imath,\tilde m,\:x, \: a_1, b_1)$ correspondingly;
\item\label{req:2} for any parameters $(a_1,b_1)\in\phi(BS(1,k))$ there exist coordinate maps $\mu_{\Gamma1}\colon \langle b_1 \rangle\to \Z$ and $\mu_{\Delta1}\colon \Z^3\to BS(1,k)$ such that the compositions $\mu_{\Gamma1}\circ\mu_{\Delta1}$ and $\mu_{\Delta1}\circ\mu_{\Gamma1}$ are definable by formulas $\theta_\Z(\tilde z,\tilde \imath, \tilde m,\:x,\:\bar r)$ and $\theta_{BS(1,k)}(\tilde z,\tilde\imath,\tilde m,\:x, \: a_1, b_1)$ for any $\bar r\in \mu_{\Delta1}^{-1}((a_1,b_1))$.
\end{enumerate}

For any parameters $(a_1,b_1)\in \phi(BS(1,k))$ we will consider two coordinate maps of the interpretation $BS(1,k)\simeq\Delta(\Z,\emptyset)$. The first one will be $\mu_\Delta$ as before~\eqref{eq:mu_Delta} ($\mu_{\Delta}\colon (z,i,m)\in \Z^3 \to (zk^i,m)\in BS(1,k)$) and the second one will be the composition of $\mu_\Delta$ and the automorphism $\lambda$ form Lemma~\ref{lemma5}, so $\mu_{\Delta1}=\lambda\circ\mu_\Delta$:
$$
\mu_{\Delta1}\colon (z,i,m)\in \Z^3 \; \longrightarrow\; a_1^{zk^i}b_1^m\in BS(1,k).
$$
And let $\mu_{\Gamma1}$ will be as before $\mu_\Gamma$~\eqref{eq:mu_Gamma}: $\mu_{\Gamma1}\colon b_1^m \in \langle b_1\rangle \to m\in \Z$. We are going to prove that $\mu_{\Gamma1}\circ \mu_{\Delta}$ and $\mu_{\Delta1}\circ \mu_{\Gamma1}$ satisfy to item~\ref{req:1}, as well as $\mu_{\Gamma1}$, $\mu_{\Delta1}$ satisfy to item~\ref{req:2} of the requirements above.

Note that the base set of the interpretation $BS(1,k)\simeq \Delta\circ\Gamma(BS(1,k),a_1,b_1)$ is 
$$
U_{\Delta\circ\Gamma}(BS(1,k), a_1,b_1)=U_\Gamma(BS(1,k),b_1)^3=\langle b_1\rangle^3,
$$
and its coordinate map is $\mu_{\Delta1}\circ\mu_{\Gamma1}=\mu_{\Delta1}\circ\mu^3_{\Gamma1}\big|_{U_{\Delta\circ\Gamma}(BS(1,k), a_1,b_1)}$:
$$
\mu_{\Delta1}\circ\mu_{\Gamma1}\colon (b_1^z,b_1^i,b_1^m)\in \langle b_1\rangle^3 \;\xrightarrow[\mu_{\Gamma1}]\; (z,i,m)\in \Z^3\;\xrightarrow[\mu_{\Delta1}]\; a_1^{zk^i}b_1^m\in BS(1,k).
$$
The graph of $\mu_{\Delta1}\circ\mu_{\Gamma1}$ is the set $\{(b_1^z,b_1^i,b_1^m,a_1^{zk^i}b_1^m)\mid z,i,m\in \Z\}$. Meanwhile, one has  $BS(1,k)\models\theta_{BS(1,k)}(\tilde z,\tilde \imath,\tilde m,\:x,\:a_1,b_1)$ if and only if $\tilde z,\tilde \imath, \tilde m\in \langle b_1\rangle$, say $\tilde z=b_1^z$, $\tilde \imath=b_1^i$, $\tilde m=b_1^m$, and $BS(1,k)\models\tau(a_1,\,b_1^i\,x\,b_1^{-m}\,b_1^{-i},\,b_1^z,\,b_1)$. But at the same time, by Corollary~\ref{cor5}, one has $BS(1,k)\models\tau(a_1,\,b_1^i\,x\,b_1^{-m}\,b_1^{-i},\,b_1^z,\,b_1)$ if and only if $b_1^i\,x\,b_1^{-m}\,b_1^{-i}=a_1^z$, i.\,e., $x=b_1^{-i}a_1^z\, b_1^i \, b_1^m=a_1^{zk^i}b_1^m$. Thus the graph of $\mu_{\Delta1}\circ\mu_{\Gamma1}$ is definable by formula 
$\theta_{BS(1,k)}(\tilde z,\tilde\imath,\tilde m,\:x, \: a_1, b_1)$.

Further, as we know, 
\begin{multline*}
\Z\models\phi_\Delta(z_a,i_a, m_a,z_b,i_b,m_b) \iff BS(1,k)\models \phi(\mu_\Delta((z_a,i_a, m_a,z_b,i_b,m_b))) \iff \\ \iff BS(1,k)\models \phi((z_ak^{i_a},m_a),(z_bk^{i_b},m_b))\iff \\ \iff (z_ak^{i_a},m_a)\in A_1, (z_bk^{i_b},m_b)\in Ab \iff z_ak^{i_a}\in U(\Z[1/k]), m_a=0, m_b=1.
\end{multline*}
Take arbitrary parameters $\bar r=(z_a,i_a, m_a,z_b,i_b,m_b)$ from $\phi_\Delta(\Z)$ and put $(a_1,b_1)=\mu_\Delta(\bar r)$, $a_1=(z_a k^{i_a},m_a)$, $b_1=(z_b k^{i_b},m_b)$. 
Since $U_{\Gamma\circ\Delta}=(U_\Gamma)_\Delta=([x,\tilde b]=e)_\Delta$, therefore, the base set $U_{\Gamma\circ\Delta}(\Z,\bar r)$ equals to $\{(z,i,m)\in\Z^3\mid [(zk^i,m),(z_bk^{i_b},m_b)]=e\}$. By Corollary~\ref{lemma3}, $C_{BS(1,k)}(b_1)=\langle (b_1)\rangle$. Since $m_b=1$ and due to~\eqref{eq:g^n}, one has $(z_bk^{i_b},m_b)^n=(z_bk^{i_b} \frac{k^{-n}-1}{k^{-1}-1},n)$ for any $n\in \Z$, so 
$$
U_{\Gamma\circ\Delta}(\Z,\bar r)=\{(z,i,m)\in\Z^3\mid zk^i=z_bk^{i_b} \frac{k^{-m}-1}{k^{-1}-1}\}.
$$
The composition $\mu_{\Gamma1}\circ\mu_{\Delta}$ is a coordinate map of the interpretation $\Z\simeq\Gamma\circ\Delta(\Z,\bar r)$:
$$
\mu_{\Gamma1}\circ\mu_{\Delta}\colon (z,i,m)\in U_{\Gamma\circ\Delta}(\Z,\bar r) \;\xrightarrow[\mu_{\Delta}]\; (zk^i,m)=b_1^m\in BS(1,k)\;\xrightarrow[\mu_{\Gamma1}]\; m\in \Z.
$$
It is clear, that the graph of $\mu_{\Gamma1}\circ\mu_{\Delta}$ is definable by formula $\theta_\Z(\tilde z,\tilde \imath, \tilde m,\:x,\:\bar r)$. 

Finally, take any tuples of parameters $(a_1,b_1)\in \phi(BS(1,k))$ and $\bar r=(z_a,i_a, m_a,z_b,i_b,m_b)\in \mu_{\Delta1}^{-1}((a_1,b_1))$. Since $\mu_{\Delta1}(\bar r)=(a_1^{z_ak^{i_a}}b_1^{m_a},a_1^{z_bk^{i_b}}b_1^{m_b})=(a_1,b_1)$, therefore, $z_ak^{i_a}=1$, $m_a=0$, $z_b=0$, $m_b=1$. Thus $U_{\Gamma\circ\Delta}(\Z,\bar r)=\{(z,i,m)\in\Z^3\mid z=0\}$. The coordinate map $\mu_{\Gamma1}\circ\mu_{\Delta1}$ of the the interpretation $\Z\simeq\Gamma\circ\Delta(\Z,\bar r)$ is defined as
$$
\mu_{\Gamma1}\circ\mu_{\Delta1}\colon (0,i,m)\in \Z^3 \;\xrightarrow[\mu_{\Delta}]\; (0,m)=b^m\in BS(1,k)\;\xrightarrow[\lambda]\; b_1^m\in BS(1,k)\;\xrightarrow[\mu_{\Gamma1}]\; m\in \Z.
$$
Thus the graph of $\mu_{\Gamma1}\circ\mu_{\Delta1}$ is definable by formula $\theta_\Z(\tilde z,\tilde \imath, \tilde m,\:x,\:\bar r)$ again, as required. 
\end{proof}

\section{Groups elementarily equivalent to $BS(1,k)$}

\subsection{Nonstandard models of arithmetic}\label{se:nonstandard_arithmetic}

Here we introduce some well-known definitions and facts on non-standard arithmetic. As a general reference we follow the book~\cite{Kaye}.

Let $L = \{+,\,\cdot\,, 0,1\}$ be the language of rings with unity $1$. By $\Z$ we denote the {\em standard arithmetic}, i.\,e., the set of integers with the standard operations from $L$. And by $\nsZ$ we denote {\em non-standard arithmetic}, i.\,e., any $L$-structure, that is elementarily equivalent to $\Z$ ($\nsZ\equiv\Z$), but not isomorphic to $\Z$. 

Let $\tilde 0$ and $\tilde 1$ be  interpretations of $0$ and $1$ in $\nsZ$. One can identify any integer number $m \in \Z$ with a non-standard number $\tilde m = \tilde 1 + \ldots + \tilde 1$ (or $\tilde m = -(\tilde 1 + \ldots  +\tilde 1$), if $m<0$), which is the sum of $m$ non-standard units $\tilde 1$ in $\nsZ$. 

The map $\lambda\colon m \to \tilde m$ gives an elementary embedding $\Z \to \nsZ$, i.\,e., $\lambda(\Z)$ is an elementary substructure of $\nsZ$, so  $\Z$ is a prime model of the  theory $\Th(\Z)$. In the sequel we always identify $\Z$ with its image $\lambda(\Z)$ via $\lambda$ and call elements of  $\lambda(\Z)$ the {\em standard integers}. We will write them as $0, 1, m, k,\ldots$, omitting tilde. It follows from Peano induction axiom that $\Z$ is not definable in $\nsZ$ (even with parameters from $\nsZ$). 

By Lagrange's four-square theorem, the set of positive integers $\Z^+$, hence the standard linear ordering $<$, is definable in $\Z$. The same formulas define positive non-standard integers $\nsZ^+$ and a linear ordering  $<$  in $\nsZ$, such that $\nsZ=-\nsZ^+\sqcup\{0\}\sqcup\nsZ^+$, where $-\nsZ^+ = (-1)\cdot \nsZ^+$.

There exists a formula $\mathbf{p}(x,y,z)$ in $L$ that defines the standard exponentiation $x=y^z$ in $\Z$ (for $z \geq 0$). The same formula defines an exponentiation $\tilde x=\tilde y^{\tilde z}$ in $\nsZ$ (where $\tilde z \geq 0$ and $\tilde x^0=1$). In particular, for any $k\in \N$ and $i\in \nsZ^+$ there exists a unique element $k^i\in \nsZ^+$, such that the following classical identities hold:
\begin{equation}\label{eq:k}
k^{i_1}\cdot k^{i_2}=k^{i_1+i_2},\quad (k^{i_1})^{i_2}=k^{i_1\cdot i_2}, \quad i_1,i_2\in\nsZ^+, 
\end{equation}
and $k^{i_1}=k^{i_2}$ implies that $i_1=i_2$. If $i\in \Z^+\subset \nsZ^+$, then $k^i$ is the usual product of $i$ factors, each equal to $k$. Thus, the set $S=\{k^i,i\in \nsZ^+\}\cup\{ 1\}$ is a multiplicative subset in $\nsZ$. 

Denote by $\nsZ[1/k^{\nsZ}]$ the ring of fractions $S^{-1}\nsZ$. Then the ring $\nsZ[1/k^{\nsZ}]$ has no zero-divisors,  $\nsZ$ embeds into $\nsZ[1/k^{\nsZ}]$, and $\nsZ[1/k^{\nsZ}]$ has characteristic zero. Further,  $\nsZ[1/k^{\nsZ}]$ is an  associative commutative unitary ring  generated by $\nsZ$ and the set of elements $S^{-1}=\{1/k^i, i\in \nsZ^+\}$ (we also write $k^{-i}=1/k^i$) with relations $k^{i_1}\cdot k^{i_2}=k^{i_1+i_2}$ for all $i_1,i_2\in \nsZ$. Any element $x\in \nsZ[1/k^{\nsZ}]$ can be written as  $x=z\cdot k^i$, where $z,i\in \nsZ$. Furthermore, for any elements $x\in S\cup S^{-1}$ and $i\in \nsZ$ one can define $x^i\in S\cup S^{-1}$, such that identities~\eqref{eq:k} hold for any $i_1,i_2\in \nsZ$.

\subsection{Non-standard models of $BS(1,k)$}
\label{subsec:5.2}

In Section~\ref{subsec:4.2} we constructed an interpretation  $BS(1,k) \simeq \Delta(\Z,\emptyset) $ which is a part of bi-interpretation of $BS(1,k)$ and $\Z$. In the sequel, we omit $\emptyset$ from the notation and write $\Delta(\Z)$.

Recall that $\Delta$ interprets $BS(1,k)$ on the set  $\Z^3$ via the coordinate map
$$
\mu_\Delta\colon (z,i,m)\in\Z^3 \;\longrightarrow \; (zk^{i},m)\in BS(1,k).
$$
Therefore, the  equivalence $\sim_\Delta$ on $\Z^3$ is defined by the rule
$$
(z_1,i_1,m_1)\sim_\Delta(z_2,i_2,m_2) \iff z_1k^{i_1}=z_2k^{i_2} \; \wedge \; m_1=m_2,
$$
while the multiplication $\odot$ and inversion $^{-1}$ on $\Z^3$ are defined by
$$
(z_1,i_1,m_1)\odot (z_2,i_2,m_2) = (z_3,i_3,m_3) \Longleftrightarrow m_3= m_1+m_2 \; \wedge \; z_3k^{i_3} = z_1k^{i_1}+z_2k^{i_2-m_1}
$$
and
$$
(z,i,m)^{-1} = (-z,i+m,-m)
$$
that correspond via $\mu_\Delta$ to the multiplication and inversion  in $BS(1,k)$. 

For every ring $\nsZ \equiv \Z$ the interpretation $\Delta$ gives a non-standard Baumslag\,--\,Solitar group $\Delta(\nsZ)$ interpreted on $\nsZ^3$ by the formulas above for $\sim_\Delta, \odot$ and $^{-1}$. 

Now, we describe the group in more algebraic terms. 
Let $k>1$ be an integer. Take the ring $\nsZ$ and construct the ring of fractions $\nsZ[1/k^{\nsZ}]$  as above. Define a metabelian group $BS(1,k,\nsZ)$ as a semidirect product 
$$
BS(1,k,\nsZ) = \nsZ[1/k^{\nsZ}]\rtimes \nsZ,
$$
where $\nsZ[1/k^{\nsZ}]$ and $\nsZ$ are viewed as the additive groups of the corresponding rings, and  $\nsZ$ acts on $\nsZ[1/k^{\nsZ}]$ via a homomorphism  $\varphi\colon \nsZ \to \Aut(\nsZ[1/k^{\nsZ}])$, where  for $m \in \nsZ$ the automorphism $\varphi(m) \in \Aut(\nsZ[1/k^{\nsZ}])$ is defined as
$$
\varphi(m)\colon x \to  x\cdot k^{-m}, \ \ x\in\nsZ[1/k^{\nsZ}].
$$

The following result is straightforward, so we omit the proof. 
\begin{lemma} Let $k>1$ and $\nsZ \equiv \Z$. Then the groups $\Delta(\nsZ)$ and $BS(1,k,\nsZ)$ are isomorphic. Moreover, a map 
$$
(z,i,m)\in\nsZ^3 \;\longrightarrow \; (zk^{i},m)\in BS(1,k,\nsZ)
,$$
gives rise to a group isomorphism $\Delta(\nsZ) \to BS(1,k,\nsZ)$.
\end{lemma}

\begin{remark}
By Theorem~\ref{th:equiv2} for any other (regular or absolute) interpretation $\Delta_1$ of $BS(1,k)$ in $\Z$ the groups $\Delta(\nsZ)$ and $\Delta_1(\nsZ)$ are isomorphic  for any $\nsZ \equiv \Z$.  Hence all of them are isomorphic to $BS(1,k,\nsZ)$. We refer to  the groups  $BS(1,k,\nsZ)$ as non-standard models of $BS(1,k)$.
\end{remark}

As we have mentioned in the Introduction the groups $BS(1,k,\nsZ)$ are interesting in their own right from algebraic, geometric, and algorithmic view-points. Here, we describe one of their interesting properties that concerns $\nsZ$-exponentiation. 

\begin{theorem} \label{th:exponent}
    For every $k>1$ and $\nsZ \equiv \Z$ the group $BS(1,k,\nsZ)$ is a $\nsZ$-group. 
    
\end{theorem}
\begin{proof}
    We mentioned in Section \ref{subsec:4.1} that  in the group $BS(1,k)$ viewed as $\Z[1/k]\rtimes \Z$, for an element $g = (y,m) \in \Z[1/k]\rtimes \Z$ and for $n \in \Z$ one has 
   \begin{equation} \label{eq:exponents}       
  g^n=(y,m)^n=
  \begin{cases}
    (y\frac{k^{-mn}-1}{k^{-m}-1},mn) & \text{if $m \neq 0$} \\
    (ny,0) & \text{if $m = 0$.}
  \end{cases}
  \end{equation}
This defines $\Z$-exponentiation in the group $BS(1,k)$. It follows that there is a formula $Exp(\bar u,\bar v,t)$, where $\bar u = (u_1,u_2,u_3), \bar v = (v_1,v_2,v_3)$, in the language of rings that defines in $\Z$ the standard exponentiation in the interpretation $\Delta(\Z)$. Namely, for any $z_1,i_1,m_1,z_2,i_2,m_2,n \in \Z$ the following equivalence holds
$$
\Z \models Exp(z_1,i_1,m_1,z_2,i_2,m_2,n) \Longleftrightarrow
\mu_\Delta(z_1,i_1,m_1)^n = \mu_\Delta(z_2,i_2,m_2).
$$
Since $\nsZ \equiv \Z$ the formula $Exp(\bar u,\bar v,t)$ in $\nsZ$ defines a $\nsZ$-exponentiation in the interpretation $\Delta(\nsZ)$ which in the group $BS(1,k,\nsZ)$ takes the form of (\ref{eq:exponents}), where 
$y\in \nsZ[1/k^{\nsZ}],m,n \in \nsZ$. Furthermore, this exponentiation satisfies the axioms of $\nsZ$-groups (1)--(4) (see Introduction), since these axioms can also be written by formulas in $\Z$ (using the formula $Exp(\bar u, \bar v, t)$). This shows that $BS(1,k,\nsZ)$ is a $\nsZ$-group.
\end{proof}


\subsection{First-order classification}
\label{subsec:5.3}

The following result describes all groups elementarily equivalent to $BS(1,k)$.

\begin{theorem} \label{th:FO classification}
Let $k$ be an integer, $k>1$. Then for any ring $\nsZ\equiv \Z$ one has 
$$
BS(1,k)\equiv BS(1,k,\nsZ)
$$ 
and, conversely, every group  elementarily equivalent to $BS(1,k)$ has the form $BS(1,k,\nsZ)$ for some $\nsZ\equiv \Z$. Moreover, this description of models of $\Th(BS(1,k))$ is a bijection: for any $\nsZ_1\equiv\Z\equiv\nsZ_2$ one has 
$$
BS(1,k,\nsZ_1) \simeq BS(1,k,\nsZ_2) \Longleftrightarrow \nsZ_1\simeq \nsZ_2.
$$
\end{theorem}

\begin{proof}
It follows from Theorems~\ref{th:equiv1} and~\ref{th:equiv2}.
\end{proof}



\bigskip
\bigskip

\noindent \textsf{\textbf{Authors:}}

\medskip

\noindent \textsf{\textbf{Evelina Daniyarova}}

Sobolev Institute of Mathematics SB RAS, Pevtsova, 13, Omsk, Russia, 644099

E-mail: \texttt{evelina.daniyarova@gmail.com}

\bigskip

\noindent \textsf{\textbf{Alexei Myasnikov}}

Schaefer School of Engineering \& Science, Department of Mathematic Sciences, Stevens Institute of Technology, Castle Point on Hudson, Hoboken NJ 07030-5991, USA

E-mail: \texttt{amiasnikov@gmail.com}

\bigskip

\end{document}